\documentclass[11pt]{article}
\input epsf.tex
\usepackage{latexsym}
\usepackage{amsmath,amsthm,amsfonts,amssymb}
\usepackage{epsfig}
\usepackage{graphicx}

\newtheorem{pro}{Proposition}[section]
\newtheorem{thm}[pro]{Theorem}
\newtheorem{lem}[pro]{Lemma}

\def\bb{{\bar b}}
\def\bx{{\bar x}}
\def\by{{\bar y}}

\def\A{{\mathbb A}}
\def\C{{\mathbb C}}
\def\mE{{\mathbb E}}
\def\mF{{\mathbb F}}
\def\H{{\mathbb H}}

\def\R{{\mathbb R}}
\def\Z{{\mathbb Z}}

\def\sech{{\textnormal{sech}}}

\def\arccosh{{\textnormal{arccosh}}}

\def\chix{{\raise.5ex\hbox{$\chi$}}}

\def\sF{{\mathcal F}}
\def\sP{{\mathcal P}}

\def\teta{{\tilde \eta}}

\def\tphi{{\tilde \phi}}

\def\tH{{\tilde H}}

\def\tG{{\tilde G}}
\def\tP{{\tilde P}}
\def\tpartial{{\tilde \partial}}

\def\tS{{\tilde S}}
\def\tgamma{{\tilde \gamma}}
\def\tsigma{{\tilde \sigma}}

\def\tv{{\tilde v}}

\title{A Generalization of the Prime Geodesic Theorem to Counting Conjugacy Classes of
Free Subgroups}

\author{Lewis Bowen \footnote{Research supported in part by a Max Zorn Postdoctoral
Fellowship}}

\begin{document}

\maketitle
\begin{abstract}
The classical prime geodesic theorem (PGT) gives an asymptotic formula (as $x$ tends to
infinity) for the number of closed geodesics with length at most $x$ on a hyperbolic
manifold $M$. Closed geodesics correspond to conjugacy classes of $\pi_1(M)=\Gamma$ where
$\Gamma$ is a lattice in $G=SO(n,1)$. The theorem can be rephrased in the following
format. Let $X(\Z,\Gamma)$ be the space of representations of $\Z$ into $\Gamma$ modulo
conjugation by $\Gamma$. $X(\Z,G)$ is defined similarly. Let $\pi: X(\Z,\Gamma)\to
X(\Z,G)$ be the projection map. The PGT provides a volume form $vol$ on $X(\Z,G)$ such
that for sequences of subsets $\{B_t\}$, $B_t \subset X(\Z,G)$ satisfying certain
explicit hypotheses, $|\pi^{-1}(B_t)|$ is asymptotic to $vol(B_t)$.

We prove a statement having a similar format in which $\Z$ is replaced by a free group of
finite rank under the additional hypothesis that $n=2$ or $3$.
\end{abstract}

%Title: A Generalization of the Prime Geodesic Theorem to Counting $\pi_1$-injective
%immersions of surfaces with boundary.

%Abstract:
%The classical prime geodesic theorem (PGT) gives an asymptotic formula (as $x$ tends to
%infinity) for the number of closed geodesics with length at most $x$ on a hyperbolic
%manifold $M$. This can be reformulated in the following way. A closed geodesic of M
%corresponds to a homotopy class of closed curves on M. This class corresponds to an
%immersion of a circle f: S^1 -> M up to homotopy. So the PGT is asymptotic formula for
%the homotopy classes of circle immersions counted according to the length of the unique
%geodesic representative. %

%In this talk I will generalize this theorem by presenting an asymptotic formula for the
%number of homotopy classes of immersions of surfaces with nonempty boundary counted
%according to the moduli of the unique geodesic representative. All the immersions
%considered are pi_1-injective. The results have been proven for dimensions 2 and 3 only
%but I conjecture that they hold in all dimensions.%


\begin{center}
\section{Introduction}
\end{center}

Given a closed manifold $M^n$ of contant curvature $-1$, the prime
geodesic theorem states that the number of closed oriented geodesics on $S$
of length at most $L$ is asymptotic to $\frac{e^{(n-1)L}}{(n-1)L}$.
This was first proven by [Huber, 1961] (in the $n=2$ case and using
Selberg's trace formula) and subsequently generalized by Margulis in
his thesis, [Parry and Pollicott, 1983], [Knieper, 1997], [Gunesch,
 2005] and others. [Hurt, 2001] and Sharp (see [Margulis, 2004]) have
written two recent surveys.

Since $M$ has negative curvature, the set of closed geodesics is naturally in 1-1
correspondence with the set of conjugacy classes of nontrivial elements of $\pi_1(M)$.
The length of a closed geodesic corresponds to the translation length
of any element in its associated
conjugacy class. So, the theorem is equivalent to an asymptotic formula for the number of
conjugacy classes of the fundamental group counted according to translation length. We
will generalize this by providing an asymptotic formula for the number of conjugacy
classes of free subgroups counted according to moduli.

%\begin{center}
\subsection{Example 1: Pairs of Pants}
%\end{center}

Before giving a general framework to describe the results, we present an example in the
simplest nontrivial case. We assume that $M$ is a closed hyperbolic surface and we are
interested in counting conjugacy classes of subgroups $F < \pi_1(M)$ where $F$ is free of
rank 2. As a subgroup of $\pi_1(M)$, $F$ corresponds to a covering space $M_F \to M$. If
$F'$ is conjugate to $F$ then $M_{F'}$ is isometric to $M_F$. Since $F$ has rank 2, $M_F$
is either an open 3-holed sphere or an open 1-holed torus.

Suppose further that we are only interested in counting conjugacy
classes $F$ such that the covering space
$M_F$ is a 3-holed sphere. Let $Mod_{0,3}$ denote the (moduli) space of all isometry
classes of
complete hyperbolic manifolds which are topologically open 3-holed spheres. There is a
natural way to parametrize $Mod_{0,3}$; there are exactly 3 pairwise nonhomotopic simple
closed curves on any 3-holed sphere. The set of lengths of their geodesic representatives
is a parametrization of $Mod_{0,3}$. So $Mod_{0,3} = (0,\infty)^3/Sym(3)$ where $Sym(3)$,
the
symmetric group on $3$ elements acts by permuting the coordinates.

Let $r_1,r_2,r_3$ be distinct positive numbers. One consequence of the main theorem in
this paper is that if $\epsilon >0$ then the number of conjugacy classes $[F]$ such that
the isometry class of $M_F$ has coordinates $\{l_1,l_2,l_3\} \in Mod_{0,3}$ satisfying
$|l_i-r_it| < \epsilon$ is asymptotic to
\begin{eqnarray*}
\frac{64\sinh(\epsilon/2)^3 e^{(r_1+r_2+r_3)t/2} }{ vol(T_1 M)}
\end{eqnarray*}
as $t \to \infty$. Here $vol(T_1(M))=2\pi(area(M))=4\pi^2 genus(M)$.

%\begin{center}
\subsection{The Framework}
%\end{center}
Let $G=SO(n,1)$ which we identify as the group of
orientation preserving isometries of the unique complete simply connected space of
constant curvature $-1$ which we will denote by $\H^n$. Then $M$ is isometric to
$\H^n/\Gamma$ for some uniform lattice $\Gamma < G$ and $\Gamma$ is isomorphic to
$\pi_1(M)$.

Let $F$ be a free group. The representation variety, denoted by $R(F,G)$, is the set of
homomorphisms $\phi:F \to G$. $G$ acts on $R(F,G)$ by conjugation $(g\phi)(f) =
g\phi(f)g^{-1}$. The quotient $R(F,G)/G$ is called the character variety and denoted by
$X(F,G)$. $R(F,\Gamma)$ and $X(F,\Gamma)=R(F,\Gamma)/\Gamma$ are defined similarly.

When $F=\Z$ the group of integers, $X(F,\Gamma)=X(\Z,\Gamma)$ is the set of conjugacy
classes of cyclic subgroups of $\Gamma$ which is naturally in 1-1 correspondence with the
closed geodesics of $M$.

$X(\Z,G)$ is the set of conjugacy classes of cyclic subgroups of $G$. It is well known
that $X(\Z,G)$ splits into 3 pairwise disjoint subsets depending on the classification of
isometries into elliptic, parabolic and hyperbolic. Hyperbolic isometries are
characterized by the property of having exactly 2 fixed points on the boundary at
infinity of $\H^n$. If $M$ is a closed manifold then all infinite order elements of
$\Gamma$ are hyperbolic.

So let $X_h(\Z,G) \subset X(\Z,G)$ be the set of conjugacy classes of hyperbolic
elements. We identify $X_h(\Z,G)$ with $(0,\infty) \times SO(n-2)$ where the first factor
accounts for the translation length of the hyperbolic element and the second factor
accounts for the holonomy about its axis. Let $vol_X$ be the measure on $X_h(\Z,G) \cong
( 0,\infty) \times SO(n-2)$ equal to the product of Lebesque measure on $(0,\infty)$ with
normalized Haar measure on $SO(n-2)$. The inclusion $\Gamma \to G$ induces a map
$\pi:X(\Z,\Gamma) \to X(\Z,G)$.

If $B_t=[t-\epsilon,t+\epsilon] \times SO(n-2) \subset X_h(\Z,G)$ then
the prime geodesic theorem is equivalent to the following:
\begin{eqnarray*}
|\pi^{-1}(B_t)| \sim \int_{B_t} \frac{e^{(n-1)x}}{x} dvol_X(x,h).
\end{eqnarray*}
Here and everywhere else in this paper $f \sim g$ means $f(t)/g(t) \to 1$ as $t \to
\infty$.

%Note that because $\Gamma$ is discrete, if $B$ is a compact subset of $X(\Z,G)$
%then $\pi^{-1}(B)$ is finite so the equation above makes sense.

The main result of this paper generalizes the formula above. It has the following form.
Let $F$ be a free group of finite rank and $G=SO(2,1)$ or $SO(3,1)$. Let $\pi:
X(F,\Gamma) \to X(F,G)$ be the map induced by inclusion $\Gamma < G$. For $t>0$, let $B_t
\subset X(F,G)$. Then, there is an explicit volume form $vol_{F,\Gamma}$ on $X(F,G)$ such
that if
certain hypotheses on $\{B_t\}$ are satisfied then $|\pi^{-1}(B_t)|
\sim vol_X(B_t)$. This volume form depends on $\Gamma$ only through
the volume of $M$.

%\begin{center}
\subsection{The Dimension 2 Case}
%\end{center}

To be more precise, we must provide appropriate hypotheses on the sequence $\{B_t\}$ and
exhibit a formula for $vol_{F,\Gamma}$. For simplicity, let us
first consider the dimension $2$ case. We will be concerned only with the subset
$X_d(F,G) \subset X(F,G)$ of discrete faithful representations with no
cusps. The last condition means that if $\phi: F \to G$ is a
representation with conjugacy class $[\phi] \in X_d(F,G)$ then $\H^2/\phi(F)$ is a
surface without cusps; i.e. $\phi(F)$ has no parabolic elements.

So $\H^2/\phi(F)$ is a surface of genus $g$ with $n$
holes for some $g$
and $n \ge 1$. The pair $(g,n)$ is invariant under conjugation by $G$ and the character
variety $X_d(F,G)$ splits as a disjoint union of $X_{g,n}(F,G)$ where $[\phi] \in
X_{g,n}(F,G)$ iff $\H^2/\phi(F)$ has genus $g$ and $n$
holes. Each connected component of $X_{g,n}(F,G)$ is identifiable with the Teichmuller
space of genus $g$ $n$-holed surfaces.

For example, if $F$ has rank 2 then $X_d(F,G) = X_{0,3}(F,G) \cup
X_{1,1}(F,G)$. Each component of $X_{0,3}(F,G)$ is identifiable with $(0,\infty)^3$ which
is the universal covering space of $Mod_{0,3}$ $= (0,\infty)^3/Sym(3)$. In general,
$Mod_{g,n}$ is the quotient of $X_{g,n}(F,G)$ by the group of automorphisms of $F$ which
acts on representations by precomposition.

Recall that a pants decomposition of a surface $S$ is a collection $\sP=\{P_i\}$ of
nonoverlapping embedded 3-holed spheres in $S$ whose union is all of $S$. A pants
decomposition determines Fenchel-Nielsen coordinates  on the Teichmuller space of $S$
[Casson and Bleiler, 1988]. Roughly speaking, the coordinates are given by a map $FN:
Teich(g,n) \to \R^k$ where $k=6g+3n-6$. By pulling back we obtain a volume form and
distance function on $Teich(g,n)$ which can then be pulled back to $X_{g,n}(F,G)$. By
[Wolpert, 1982] this volume form coincides with the
Weil-Petersson volume form. So it is independent of the pants decomposition. We denote it
by $vol_{WP}$. The distance function depends on the pants decomposition. We denote it by
$d_\sP$. In section \ref{sec:3d} below we explain these coordinates in greater detail and
extend them to what we call the 'nondegenerate' character varieties
$X_\sP(F,Isom^+(\H^n))$ for
$n=2,3$.

By definition the convex core of a surface is the smallest convex subset that is homotopy
equivalent to the surface. If $[\phi] \in X_d(F,G)$, the convex core of $\H^2/\phi(F)$ is
a compact surface with geodesic boundary. Let $l_\partial([\phi])$ denote the total
length of this boundary.

For example, if $\H^2/\phi(F)$ is an open 3-holed sphere, then its convex core is a
compact
3-holed sphere and $l_\partial([\phi])=l_1+l_2+l_3$ where $\{l_1,l_2,l_3\} \in Mod_{0,3}$
are the coordinates for the isometry class of the surface $\H^2/\phi(F)$.

If for each $t>0$, $[\phi_t] \in X_d(F,G)$ and $\sP_t=\{P_{t,i}\}_{i=1}^{r-1}$ is a pants
decomposition of the
convex core of $\H^2/\phi_t(F)$ (where $r$ is the rank of $F$) then we call the path
$\{(\phi_t,\sP_t)\}_{t>0}$ {\bf
telescoping} if the following holds.
\begin{itemize}

\item Let $S_{t,i}:= \cup_{j<i} P_{t,j}$. Then $P_{t,i} \cap S_{t,i}$ has either $0, 1$
or $2$ components.

\item For all $i>1$, $\lim_{t \to \infty} len(\partial S_{t,i}) -
len(\partial S_{t,i-1}) = + \infty$.
\end{itemize}

%It is proven in section \ref{sec:telescope} that if $F$ is a free group of finite rank
%and $[\phi_t] \in X_d(F, Isom^+(\H^2))$ is a path such that the length of the shortest
%closed curve on $\H^2/\phi_t(F)$ tends to infinity then there exists a decomposition
%$\sP_t$ for which $(\phi_t,\sP_t)$ is telescoping. It is possible that an analogous
%result holds in all dimensions but we were unable to prove it.

As a corollary to the main result we obtain:

\begin{thm}\label{thm:2d}
Suppose $\{(\phi_t, \sP_t)\}_{t>0}$ is telescoping where $[\phi_t] \in
X_{g,n}(F,Isom^+(\H^2))$.
Then
\begin{eqnarray*}
|\pi^{-1}(N_\epsilon \, [\phi_t] )| \sim vol(T_1 M)^{1-rank(F)} 2^{-g} \int_{N_\epsilon
\, [\phi_t]} e^{l_\partial(\psi)/2}
\, dvol_{WP}(\psi).
\end{eqnarray*}
Here $\pi: X(F,\Gamma) \to X(F,G)$ is the projection map and $N_\epsilon \, [\phi_t]$
denotes the $\epsilon$-neighborhood of $[\phi_t] \in X_d(F,G)$ with respect to the metric
$d_{\sP_t}$.
\end{thm}
%The metric on $X_0(F,G)$ used above is described in section \ref{sec:distance} below. We
%conjecture that
%the theorem is true with respect to the Teichmuller metric but
%cannot prove it at this time.
{\bf Remark:} The hypotheses are not as restrictive as they might look.
In order to obtain any asymptotics of the above sort it is necessary to assume that the
length of the shortest closed geodesic of the surface $\H^2/\phi_t(F)$ tends to infinity.
For example if the shortest geodesic of $\H^2/\phi_t(F)$ is much shorter than the shortest
geodesic of $M$ then $|\pi^{-1}(N_\epsilon \, [\phi_t] )| =0$. But we will show in
section \ref{sec:telescope} that if the length of the shortest geodesic of
$\H^2/\phi_t(F)$ tends to infinity then there exists a decomposition $\sP_t$ for which
$\{(\phi_t,\sP_t)\}_{t >0}$ is telescoping. It is possible that a similar result holds in
all
dimensions but we were unable to prove it. Our main result proves the above formula for
sequences of subsets more general than
$\{N_\epsilon \, [\phi_t]\}$.

\subsubsection{Example 2: One Holed Tori}

If $M$ is a closed hyperbolic surface then how many immersed 1-holed
tori are there with geodesic boundary of length $b \in [L-\epsilon,L+\epsilon]$?

To answer this, recall that if $F$ is a free group of rank 2 then
$X_{1,1}(F,\Gamma)$ is the set of conjugacy classes of faithful
homomorphisms $\phi: F \to \Gamma$ such that the corresponding
covering space $M_F$ is a one-holed torus. Every component of
$X_{1,1}(F,\Gamma)$ is identifiable with the Teichmuller space of the
one-holed torus. The Teichmuller space orbifold-covers the moduli
space of the one-holed torus. So fix a component $X^0_{1,1}(F,\Gamma)
\subset X_{1,1}(F,\Gamma)$ and let $Z \subset
X^0_{1,1}(F,\Gamma)$ be a fundamental domain for the action of the
mapping class group on $X^0_{1,1}(F,\Gamma) = Teich_{1,1}$. So, under
the orbifold cover $Teich_{1,1}\to Mod_{1,1}$, $Z$
projects onto the moduli space $Mod_{1,1}$ in a 1-1 way everywhere
except on its boundary which we may assume is piecewise smooth.

For each
$b>0$, let $Z_b \subset Z$ be the set of isometry classes of one-holed
tori whose boundary has length $b$. If $Y$ is a subset of the real
line let $Z_Y = \cup_{y \in Y} \, Z_y$. In the language of the above
theorem, the number of immersed 1-holed tori with geodesic boundary of
length $b \in [L-\epsilon,L+\epsilon]$ equals
$|\pi^{-1}(Z_{[L-\epsilon,L+\epsilon]})|$.

We cannot directly apply the theorem above for
two reasons. First, the sets $Z_{[L-\epsilon,L+\epsilon]}$ are not in
the form of an $\epsilon$-neighborhood. Second, we do not have a
telescoping decompositions.

The first problem can be overcome by appealing to the more general
hypotheses of the main theorem (theorem \ref{thm:3d}). This basically amounts to covering
$Z_{[L-\epsilon,L+\epsilon]}$ with
$\delta$-neighborhoods for some $\delta>0$. To remove the second
problem, for $\delta>0$, let $Z^\delta_{b}$ be the set of isometry
classes of one-holed tori in $Z_b$
such that the length of the shortest closed geodesic is at least the product
$\delta b$. By the remark above, if for each $t>0$, $[\phi_t]$ is any
arbitrary element of $Z_{t}^\delta$ then
there exists a pants decomposition $\sP_t$ such that the 'path' $\{(\phi_t,\sP_t)\}$
is telescoping. So we can use the sets $Z^\delta_{[L-\epsilon,L+\epsilon]}$ instead of
$Z_{[L-\epsilon,L+\epsilon]}$. The volume of $Z^\delta_b$ divided by the volume of $Z_b$
tends
to a constant (depending on $\delta$) as $b \to \infty$. This constant tends to 1 as
$\delta
\to 0$. Since we are only interested in asymptotics, a diagonal
argument implies we can apply the formula directly to $Z_b$. So we obtain
\begin{eqnarray*}
|\pi^{-1}(Z_{[L-\epsilon,L+\epsilon]})| \sim  \frac{1}{2 vol(T_1 M)}
\int_{L-\epsilon}^{L+\epsilon}
\, vol( Mod_{1,1}(b) )\, e^{b/2} \, db
\end{eqnarray*}
where $vol( Mod_{1,1}(b))$ is the Weil-Petersson volume of
$Mod_{1,1}(b)$. According to [Mirzakhani, 2006] $vol(Mod_{1,1}(b)) = b^2/24 + \pi^2/6 \sim
b^2/24$. Hence
\begin{eqnarray*}
|\pi^{-1}(Z_{[L-\epsilon,L+\epsilon]})| \sim
\frac{\sinh(\epsilon/2)L^2 e^{L/2}}{ 12vol(T_1 M)}.
\end{eqnarray*}

For comparison, recall that the number of closed oriented geodesics
with length in $[L-\epsilon,L+\epsilon]$ is asymptotic to $\frac{2\sinh(\epsilon)e^L}{L}$.

\subsubsection{Example 3: Two Holed Tori}
Let's assume that $M$ is a closed hyperbolic surface. There is a special reason why
someone might want to count immersed two-holed tori in $M$.
To explain, we introduce the 2-holed torus graph $G_M=(V,E)$ of $M$. The vertex set $V$ of
$G_M$ is the set of oriented closed geodesics of $M$. There is a directed edge from $v_1$
to $v_2$ iff there exists an immersed $\pi_1$-injective 2-holed torus $T$ such that if
$T$ is given the
orientation induced by $M$ then the boundary of $T$ is $v_1 \cup -v_2$ (where $-v_2$
means $v_2$ with the opposite orientation). A directed cycle in $G_M$ corresponds to an
immersed {\it closed} surface $S$ in $M$, obtained by gluing the 2-holed tori together
which individually represent edges in the cycle. The
fundamental group of $S$ injects into the fundamental group of
$M$.

The definition of a 2-holed torus graph can be generalized to the
case when $M$ is not a surface except that the vertices then correspond to
unoriented geodesics and the edges to (unoriented) immersed $\pi_1$-injective
2-holed tori. Cycles still correspond to immersed surfaces but these need not be
$\pi_1$-injective.

It can be shown using the methods in this paper that every vertex of $G_M$ has infinite
valence. If $M$ contains an immersed $\pi_1$-injective closed surface, then $G_M$
necessarily contains directed cycles. It would be interesting to have a more
ergodic-theoretic proof that $G_M$ contains directed cycles when $M$ is closed surface.
Note that path components of $G_M$ correspond to homology classes of $M$.

For $L, \epsilon > 0$ let $G_{L,\epsilon}$ be the subgraph of $G$ induced by the set of
vertices $v$ whose
underlying oriented geodesic has length in the interval
$[L-\epsilon,L+\epsilon]$.

%Next we obtain the asymptotic number of vertices and edges for $G_{L,\epsilon}$ and
therefore the
%average degree of a vertex.

The classical prime geodesic theorem implies that the number of vertices of
$G_{L,\epsilon}$ is asymptotic to $2\sinh(\epsilon)e^L/L$. To obtain
asymptotics for the number of edges of $G_{L,\epsilon}$ we proceed as
in the previous example.

So fix a component $X^0_{1,2}(F,\Gamma)
\subset X_{1,2}(F,\Gamma)$ and let $Z \subset
X^0_{1,2}(F,\Gamma)$ be a fundamental domain for the action of the
mapping class group on $X^0_{1,2}(F,\Gamma) = Teich_{1,2}$.

For each
$(b_1,b_2)>0$, let $Z(b_1,b_2) \subset Z$ be the set of isometry classes of two-holed
tori with geodesic boundary components of length $b_1$ and $b_2$. If
$Y \subset \R^2$ let $Z_Y = \cup_{y \in Y} \, Z(y)$. In the language of the above
theorem, the number of immersed 2-holed tori with geodesic boundary lengths $(b_1,b_2)
\in [L-\epsilon,L+\epsilon]^2$ equals
$|\pi^{-1}(Z([L-\epsilon,L+\epsilon]^2))|$.

As before, we cannot directly apply the theorem above because the sets we are interested
in are not $\epsilon$
neighborhoods and do not necessarily have telescoping pants
decompositions. But the same tricks used in the previous theorem apply
here as well. Therefore
\begin{eqnarray*}
|\pi^{-1}(Z_{[L-\epsilon,L+\epsilon]^2})| \sim  \frac{1}{2 vol(T_1 M)^2}
\int_{L-\epsilon}^{L+\epsilon}\int_{L-\epsilon}^{L+\epsilon}
\, vol( Mod_{1,2}(b_1,b_2) )\, e^{(b_1+b_2)/2} \, db_1db_2
\end{eqnarray*}
where $vol( Mod_{1,2}(b_1,b_2))$ is the Weil-Petersson volume of
$Mod_{1,2}(b_1,b_2)$. According to [Mirzakhani, 2006]
$vol(Mod_{1,2}(b_1,b_2))$ is a polynomial of degree 4 in $(b_1,
b_2)$. The leading coefficients of this polynomial are given by a
recursive formula which seems feasible to compute but we did not do
it. So all we can say is that there is some constant $C>0$ such that
\begin{eqnarray*}
|\pi^{-1}(Z([L-\epsilon,L+\epsilon]^2))|\sim \frac{C \sinh^2(\epsilon/2) L^4
e^L}{vol(T_1(M))^2}.
\end{eqnarray*}
Therefore, the average degree of a vertex in the graph
$G_{L,\epsilon}$ is
\begin{eqnarray*}
\frac{CL^5\sinh^2(\epsilon/2)}{\sinh(\epsilon)vol(T_1(M))^2}.
\end{eqnarray*}
It would be interesting to know the degree sequence of this graph.

\subsection{Proof Sketch of Theorem \ref{thm:2d}}

A discrete faithful homomorphism $\phi: F \to \Gamma$ determines a subgroup $\phi(F) <
\Gamma$. Let $\tphi: \tS \to M$ be the associated covering space. $\tS$ is naturally
endowed with a hyperbolic metric so that $\tphi$ is a local isometry. Let $S$ be the
convex core of $\tS$. The map $\tphi: S \to M$ is a locally isometric immersion.

If $\phi'$ is $\Gamma$-conjugate to $\phi$, then the induced immersion $\tphi':S' \to M$
is related to $\tphi:S \to M$ by an isometry $\Psi: S \to S'$ in the sense that $\tphi'
\circ \Psi = \tphi$. Thus up to this natural equivalence relation, a character $[\phi]
\in X(F,\Gamma)$
determines a locally isometric immersion $\tphi:S \to M$. The converse is also true. So
the problem of counting characters is equivalent to the problem of counting locally
isometric immersions.

By decomposing $S$ into pants we see that it suffices to count locally isometric
immersions of a pair of pants $P \to M$ which satisfy various geometric and boundary
constraints. To be more precise, we need to answer questions of the following type.

Given a closed geodesic $\gamma$ in $M$, what is the number of immersed pants $P$ in $M$
with one boundary component of $P$ equal to $\gamma$ and such that the other boundary
components have lengths in prespecified intervals? Because we also care about twist
parameters, we need to be able to count the number of such pants with additional
restrictions on the relative positions of the other boundary components with respect to
$\gamma$.

Another type of question we need to answer is as follows. Given two closed geodesics
$\gamma_1, \gamma_2 \in M$, what is the number of immersed pants $P$ in $M$ with two of
its boundary components equal to $\gamma_1$ and $\gamma_2$ and such that the other
boundary component has length in a prespecified interval? Twist parameters come into play
here as well.

To answer both questions observe that a hyperbolic structure on a pair of pants is
determined by the lengths of its 3 boundary components. But the hyperbolic structure is
also determined by the lengths of just 2 boundary components if the length of the
shortest arc between those two is already known. The shortest arc is perpendicular to
both components. So rather than counting pants immersions directly, it suffices to count
"pairs of eyeglasses"; a pair of eyeglasses is a
triple $(\sigma_1,\sigma_2,\gamma)$ such that $\sigma_1, \sigma_2$ are closed geodesics in
$M$ and $\gamma$ is a segment that is perpendicular to both $\sigma_1$ and $\sigma_2$ at
its endpoints. Such a triple uniquely determines a locally isometric immersion of a pair
of pants into $M$. (Indeed, the pair of pants is homotopy equivalent to a regular
neighborhood of $\sigma_1 \cup \sigma_2 \cup \gamma$).

Both questions can be further reduced (by a change of variables), to counting the number
of perpendiculars between two segments. That is, we now suppose $\sigma_1, \sigma_2$ are
oriented geodesic segments in $M$ and ask how many perpendiculars are there between
$\sigma_1$ and $\sigma_2$ with length in some prespecified interval. It seems possible
that this question has appeared elsewhere but we were unable to find it in the literature.

We describe the pertinent result here in a little more detail since it may be of
independent interest. We allow $\sigma_1$ and $\sigma_2$ to depend on a parameter we
denote by $L$.
For $i=1,2$ let $T_1^\perp \sigma_i$ be the set of unit tangents vectors perpendicular to
$\sigma_i$. Any segment $\gamma$ perpendicular to $\sigma_1$ and $\sigma_2$ determines
vectors $v_i(\gamma) \in T_1^\perp(\sigma_i)$ tangent to $\gamma$ at its endpoints.

Let $\chi=\chi_{\sigma_1,\sigma_2}$ be the measure on $T_1^\perp(\sigma_1)\times
T_1^\perp(\sigma_2)\times [0,\infty)$ defined by setting $\chi(E)$ equal to the number of
perpendicular segments $\gamma$ with $(v_1(\gamma),v_2(\gamma),length(\gamma)) \in E$.
Let $\chi'$ be the measure on $T_1^\perp(\sigma_1)\times T_1^\perp(\sigma_2)\times
[0,\infty)$ with density
\begin{eqnarray*}
d\chi' &=&  \frac{e^{\Re(L)}}{2 vol(T_1 M)}  dvol_{T_1^\perp(\sigma_1)} \,
dvol_{T_1^\perp(\sigma_2)} \, dL.
\end{eqnarray*}

$vol_{T_1^\perp(\sigma_i)}$ is the obvious measure on $T_1^\perp(\sigma_i) = \sigma_i
\times S^0$ with total mass $2length(\sigma_i)$ and $vol(T_1 M)=(2\pi) area(M)$. Let
$\sigma^+_i =\sigma_i \times \{+1\} \subset T_1^\perp(\sigma_i)$. We will prove that if
$\epsilon>0$ is fixed and $length(\sigma_i(L))e^{L/2}$ tends to infinity with $L$ for
both $i=1,2$ then
\begin{eqnarray*}
\chi\left(\sigma^+_1 \times \sigma^+_2 \times [L-\epsilon,L+\epsilon]\right) \sim
\chi'\left(\sigma^+_1 \times \sigma^+_2 \times [L-\epsilon,L+\epsilon]\right).
\end{eqnarray*}
Here we use the notation $f(L) \sim g(L)$ to mean $\lim_{L \to \infty} \frac{f(L)}{g(L)}
= 1$. An analogous statement holds in dimension 3. See theorem \ref{thm:perp} for a
precise statement.

In summary, we first prove an asymptotic formula for the number of perpendiculars between
two segments. Next we introduce detailed notation for describing the geometry of
characters in $X(\pi_1(P),\Gamma)$. From this it is easy to answer questions about the
number of pants immersions with two boundary components fixed. A longer change of
variables argument is necessary to handle the case in which only one boundary component is
fixed. The main theorem follows from an inductive argument on the pants decomposition.

\subsection{Notation and Results for Dimensions 2 and 3}\label{sec:3d}

We will use the following notation throughout the paper. Identify $M$ with the quotient
$\H^n/\Gamma$ where $\Gamma \equiv \pi_1(M)$ is a discrete
subgroup of $G=Isom(\H^n)$. If $n=2$ set $G:=PSL_2(\R)$. If $n=3$ set $G:=PSL_2(\C)$. Let
$\A:=\C/<2\pi i>$ and let $\mF=\R$ or $\A$ depending on whether $n=2$ or $3$.
\subsubsection{ Lengths and Widths}

The length $len(g) \in \mF$ of
a hyperbolic element $g \in G$ is defined by
\begin{eqnarray*}
\cosh(len(g)/2) &=& \pm trace(g)/2\\
\Re(len(g)) &\ge& 0.
\end{eqnarray*}
Alternatively, $\Re(len(g))$ is the minimum of $d(x,gx)$ over all $x \in \H^n$ where
$d(\cdot,\cdot)$ is the distance function. $\Im(len(g))$ measures the amount of turning
along the axis of $g$.

A {\bf double cross} is a triplet $(\sigma_1, \sigma_2; \eta)$ where $\sigma_1, \sigma_2,
\eta$
are oriented geodesics in $\H^3$ such that $\eta$ intersects $\sigma_1$ and $\sigma_2$
orthogonally \cite{Fen1989}. Let $g \in G$ be such that $g(\eta)=\eta$ and
$g(\sigma_1)=\sigma_2$. Then
the {\bf width} of $(\sigma_1,\sigma_2;\eta)$ is defined to be $\pm len(g)$ with the sign
determined by:
\begin{itemize}
\item if $\sigma_1 \cap \sigma_2 = \emptyset$ then the sign is positive iff $\eta$ is
directed from $\sigma_1$ to $\sigma_2$,
\item otherwise, the sign is positive iff $(\sigma_1, \sigma_2, \eta)$ is a positively
oriented frame as determined by the right-hand rule.
\end{itemize}

%Let $S$ be a surface with boundary. The representation variety $R(\pi_1(S),G)$ is
%the space of homomorphisms $\phi:\pi_1(S) \to G$. $G$ acts on $R$ by
%%conjugation: $(g\phi)(f)=g\phi(f)g^{-1}$. The character variety $X(\pi_1(S),G)$
%equals $R$ modulo this action.

\subsubsection{Fenchel-Nielsen Coordinates and a Volume Form on the Character Variety}

Let $S$ be a compact surface with boundary. We will have use for coordinates on
$X(\pi_1(S),G)$ that are similar to the Fenchel-Nielsen
coordinates [Casson and Bleiler, 1988] on Teichmuller space. To describe these we
introduce the notion of an oriented framed pants
decomposition.

A pants decomposition of $S$ is a collection $\sP=\{P_i\}_{i=1}^{r-1}$ of pairwise
nonoverlapping embedded pants $P_i$ in $S$ whose union is all of $S$. $r$ is the rank of
the free group $\pi_1(S)$. We say that a representation $\phi: \pi_1(S) \to Isom^+(\H^n)$
is {\bf nondegenerate} with respect to $\sP$ iff
\begin{itemize}
\item for every $\gamma \in Curve(\sP)$, $\phi([\gamma])$ is a hyperbolic isometry and
\item if $\gamma_1, \gamma_2$ are boundary components of some pair of pants $P \in \sP$
and $[\gamma_1], [\gamma_2] \in \pi_1(S)$ represent them and freely generate $\pi_1(P) <
\pi_1(S)$ then the axes of $\phi([\gamma_1])$ and $\phi([\gamma_2])$ do not intersect at
infinity.
\end{itemize}
This condition is invariant under conjugation by $G$, so it makes sense to say that a
character $[\phi] \in X(\pi_1(S),G)$ is nondegenerate wrt $\sP$ iff $\phi$ is
nondegenerate wrt $\sP$. Let $X_\sP(\pi_1(S),G)$ denote $\sP$-nondegerate characters.
Note that $[\phi]$ is nondegenerate wrt $\sP$ iff $\phi$ restricted to $P$ is
nondegenerate wrt $\{P\}$ for all pants $P \in \sP$.

If $\phi$ is faithful and $\phi(\pi_1(S))$ is discrete and contains no parabolic elements
then $\phi$ is nondegenerate with respect to any pants decomposition. To see this, note
that the first condition above is automatically satisfied. The second one follows from
the fact that if $g_1, g_2 \in Isom(\H^n)$ are two hyperbolic isometries that share a
fixed point then the group generated by $g_1$ and $g_2$ is either cyclic, nondiscrete, or
it contains an elliptic or parabolic isometry.

For any $i$, any boundary component of $P_i$ is called a {\bf curve} of $\sP$. It is a
boundary curve if it is also a
boundary component of $S$. Otherwise it is an interior curve. Let $Int(\sP)$ denote the
set
of interior curves and $Curve(\sP)$ the set of all curves. An orientation on $\sP$ is a
choice of orientation for each curve of $\sP$. Fixing an orientation, a framing of $\sP$
is a function
$f: Int(\sP) \to Curve(\sP)\times Curve(\sP)$ such that if $P_1$ is the pair of pants to
left of $\gamma$, $P_2$ is the pair of pants to the right of $\gamma$ and
$f(\gamma)=(\gamma_1,\gamma_2)$ then $\gamma_i \in P_i$ for each $i$.

Fix an oriented framed pants decomposition $(\sP, f)$. We will define an $\A$-valued
twist parameter for each interior curve $\gamma$ with respect to a nondegenerate
representation $\phi:
\pi_1(S) \to G$ generalizing the usual twist parameter.

Implicit in the definition of $\pi_1(S)$ is a basepoint *. After conjugating if
necessary, we may assume that * is on the curve $\gamma$. For $i=1,2$, let $\gamma'_i$ be
an oriented closed path based at * freely homotopic to $\gamma_i$ in $P_i$. Let
$[\gamma'_i] \in \pi_1(S,*)$ represent $\gamma'_i$. By choosing the paths $\gamma'_i$
appropriately we may assume that $[\gamma'_i]$ and $[\gamma]$ freely generate
$\pi_1(P_i,*)$.

If $\phi([\gamma'_i])$ and $\phi([\gamma])$ are hyperbolic elements of $Isom(\H^n)$ then
let $\tgamma_i$ be the axis of
$\phi([\gamma'_i])$ and $\tgamma$ the axis of $\phi([\gamma])$.

If $\tgamma$ and $\tgamma_i$ do not share a point on the boundary $\partial_\infty \H^n$
let $\tsigma_i$ be the unique geodesic perpendicular to both $\tgamma$ and $\tgamma_i$.
We give $\tsigma_i$ a canonical orientation by requiring that the width $w_i$ of
$(\tgamma_i,\tgamma;\tsigma_i)$ has nonnegative real part and in case it is pure
imaginary, then it has a representative in $[0,\pi]i\subset \A$.

Finally, define the twist parameter
$twist_\phi(\gamma)=width(\tsigma_1,\tsigma_2;\tgamma)$. It is undefined in degenerate
cases but this will not
concern us in the sequel.

For the representation $\phi: \pi_1(S)\to G$ define $len_\phi:
Curve(\sP) \to \mF$ by $len_\phi(\gamma) = len(\phi \, [\gamma] )$ where $[\gamma]\in
\pi_1(S)$ is any element representing $\gamma$. The two functions $len_\phi$ and
$twist_\phi$ are invariant under the action of $G$ by conjugation.
Therefore, they give rise to functions $len_{[\phi]}, twist_{[\phi]}$ that depend only on
the conjugacy class $[\phi] \in X_\sP(\pi_1(S),G)$. Define
\begin{eqnarray*}
len_\phi(\partial S) = \sum_{\gamma \in Curve(\sP)-Int(\sP)} \, len_\phi(\gamma).
\end{eqnarray*}

By definition, $X_\sP(\pi_1(S),G)$ is the subset of the character variety on which the
length and twist functions are well-defined. $X_\sP(\pi_1(S),G)$ is identified with a
subset of $\mF^{Curve(\sP)\sqcup Int(\sP)}$ through
$[\phi] \to len_{[\phi]} \sqcup twist_{[\phi]}$. Let $vol_{\sP}(\cdot)$ be the pullback
of the volume form on $\mF^{Curve(\sP)\sqcup Int(\sP)}$.

%Let $X_0=X_0(\pi_1(S),G) \subset X(\pi_1(S),G)$ be the subset on which both functions are
%well-defined (e.g. there are no degeneracies as in the definition of the twist
%parameter). In the $2$-dimensional case we also require that $[\phi] \in
%X_0(\pi_1(S),G)$ is $\pi_1$-injective and that the twist parameters are real.

In the 2-dim case $X_d(\pi_1(S),G) \subset X_\sP(\pi_1(S),G)$ and the restriction of
$vol_\sP$ to $X_d(\pi_1(S),G)$ is the Weil-Petersson volume. We do not know whether a
similar statement holds in dimension 3 or in any higher dimension.

\subsubsection{ Telescoping Paths in the Character Variety}

Let $\sP=\{P_{i} \}_{i=1}^{r-1}$ be an ordered pants decomposition of $S$. For $t >0$
let $[\phi_t] \in X_d(\pi_1(S),G)$.  We say that $\{(\phi_t,\sP)\}_{t>0}$ is {\bf telescoping} if
the following
holds.
\begin{itemize}

\item Let $S_i$ be the surface $\cup_{j<i} P_{j}$. Then $P_{i} \cap S_{i}$ has either $0,
1$ or $2$ components.

\item For all $i>1$, $\lim_{t \to \infty} len_{\phi_t}(\partial S_{i}) -
len_{\phi_t}(\partial S_{i-1}) = \infty$.
\end{itemize}

%We will prove in section \ref{sec:telescope} that if $F$ is a free group of finite rank
%and $[\phi_t] \in X_0(F, Isom^+(\H^2))$ is a sequence of discrete and faithful characters
%such that the length of the shortest closed curve on $\H^2/\phi_t(F)$ tends to infinity
%then there exists a decomposition $\sP$ for which $(\phi_t,\sP)$ is telescoping. It is
%possible that this result holds in all dimensions but we were unable to prove it.

\subsubsection{Neighborhoods}
Let $\sP=\{P_{i} \}_{i=1}^{r-1}$ be an ordered pants decomposition of $S$. Let $S_i$ be
as in the previous subsection.

Given a set $E \subset X_\sP(\pi_1(S),G)$ and an $\epsilon >0$ there are two different
$\epsilon$-neighborhoods of the boundary $\partial E = {\bar E} - int (E)$ that we will
consider. The first is $N^\sP_1(\epsilon, \partial E)$ equal to the set of all $\psi\in
X_\sP(\pi_1(S),G)$ such that there exists $\phi \in \partial E$ satisfying the following.
\begin{itemize}
\item $l_\psi(\gamma)=l_\phi(\gamma)$ for all curves $\gamma$ in the decomposition $\sP$.
\item For any $i$ if $\gamma$ is a curve in $S_{i+1}-S_i$ then
\begin{eqnarray*}
\Big| \Re \big( twist_\psi(\gamma)-twist_\phi(\gamma) \big)\Big| &\le&  \epsilon
\exp\big(\Re(len_\psi \partial S_{i}-len_\psi \partial S_{i+1}) /4\big)\\
\Big| \Im \big( twist_\psi(\gamma)-twist_\phi(\gamma) \big)\Big| &\le&  \epsilon
\exp\big(\Re(len_\psi \partial S_{i}-len_\psi \partial S_{i+1}) /4\big).
\end{eqnarray*}
\end{itemize}

The second neighborhood is $N^\sP_2(\epsilon, \partial E)$ equal to the set of all
$\psi\in X_\sP(\pi_1(S),G)$ such that there exists $\phi \in \partial E$ with
\begin{itemize}
\item $|len_\psi(\gamma)=l_\phi(\gamma)| < \epsilon$ for all curves $\gamma$ in the
decomposition $\sP$,
\item $twist_\gamma(\psi)=twist_\gamma(\phi)$ for all curves $\gamma$ in the
decomposition $\sP$.
\end{itemize}

\subsubsection{ Statement of the Main Result}

There is a natural projection map $\pi:X(\pi_1(S),\Gamma) \to X(\pi_1(S),G)$ induced by
inclusion $\Gamma < G$. Let $r=rank(\pi_1(S))$ and $g=genus(S)$. If $n=2$ then let $\nu$
be the measure on $X_\sP(\pi_1(S),G)$ with derivative
\begin{eqnarray*}
d\nu &=& vol(T_1 M)^{1-r} 2^{-g} e^{\Re[len_{\psi}(\partial S)]/2} \, dvol_{\sP}(\psi).
\end{eqnarray*}
If $n=3$, define $\nu$ by
\begin{eqnarray*}
d\nu &=& vol(T_1 M)^{1-r} (2\pi)^{-r} e^{\Re[len_{\psi}(\partial S)]} \, dvol_{\sP}(\psi).
\end{eqnarray*}

The main result is:

\begin{thm}\label{thm:3d}
Let $\sP$ be a pants decomposition of $S$ and $\{E_t\}$ a path of subsets $E_t \subset
X_\sP(\pi_1(S),G)$. If for every $[\phi_t] \in E_t$, the path $\{(\phi_t,\sP)\}_{t>0}$ is
telescoping and if
\begin{eqnarray*}
\lim_{t \to \infty} \frac{\nu(N_1(1, \partial E_t))}{\nu(E_t)} =0 & \textrm{ and } &
\lim_{\epsilon \to 0} \limsup_{L \to \infty} \frac{\nu(N_2(\epsilon, \partial
E_t))}{\nu(E_t)} =0
\end{eqnarray*}
then $|\pi^{-1}(E_t)| \sim \nu(E_t )$ as $t$ tends to infinity.
\end{thm}

{\bf Acknowledgements}
I'd like to thank Joel Hass for inspiring conversations which helped this work get
started. I'd like to thank Chris Judge and Chris Connell for conversations and
encouragement on this subject.

\section{Notations and Conventions}\label{sec:note}

Throughout the paper, $M^n$ denotes a fixed oriented closed hyperbolic
$n$-manifold and $T_1 M$
its unit tangent space. Let $FM$ be the space of positively oriented orthonormal
$n$-frames. Each $n$ frame $(v_1,..,v_n) \in T_1 M\times...\times T_1 M$ projects to the
same point $x \in M$ under the canonical projection from $T_1 M$ to $M$. The frames $FM$
form a fiber bundle over $T_1 M$ with a natural projection $\pi:FM\to T_1 M$ which
selects the first vector in the frame, that is, $\pi(v_1,..,v_n)=v_1$. $Isom^+(\H^n)$
acts on the frame bundle of $\H^n$ on the left by
\begin{eqnarray*}
G(v_1,..,v_n) = (DG_x(v_1),...,DG_x(v_n))
\end{eqnarray*}
where $x$ is the common basepoint of $v_1,...,v_n$. Fix a reference frame $f_0 \in \H^n$.
Identify the frame $Gf$ with $G$ itself. Because of this identification, the right action
of $Isom^+(\H^n)$ on $F(\H^n)$ is defined by $(Hf_0)G=HGf_0$ for any $H,G \in
Isom^+(\H^n)$.

Since $Isom^+(\H^n)=F(\H^n)$, $FM$ is identified with $\Gamma \backslash Isom^+(\H^n)$
and $M$ itself with $\Gamma \backslash \H^n$ where $\Gamma < Isom^+(\H^n)$ is a uniform
lattice.

Let $\A = \C/<2\pi i>$. Let $\mF = \R$ or $\A$ depending on whether
$n=2$ or $n=3$. For $L \in \mF$ let
\begin{displaymath}
G_L=\left[ \begin{array}{cc}
e^{L/2} & 0 \\
0 & e^{-L/2}
\end{array} \right].
\end{displaymath}
If $n=2$, $G_L \in PSL_2(\R)=Isom^+(\H^2)$. If $n=3$, $G_L \in PSL_2(\C)=Isom^+(\H^3)$.
If $L>0$ is real and positive, the map on $FM$ defined by $f\to fG_L$ is the time $L$ map
of the frame flow. Its projection to the unit tangent bundle $T_1 M$ is the time $L$ map
of the geodesic flow.

For $L, \epsilon  \in \mF$, define $B(L,\epsilon)$ to be the set of $x \in \mF$ such that
$|\Re(x)-\Re(L)| \le |\Re(\epsilon)|$ and, in case $n=3$, $\min_{n \in \Z} \,
|\Im(x)-\Im(L) + 2\pi n| \le |\Im (\epsilon)|$ where $|\Im(\epsilon)|$ denotes the number
in $[0,\pi]$
equal to the distance from $\Im(\epsilon)$ to $0$ mod $2\pi$.

%\begin{displaymath}
%B(L,\epsilon)= \left\{\begin{array}{ll}
% \{ x \in \R ,\, |L-x| < |\epsilon| \}& \textrm{ , if }
%n=2\\
%\{x+iy \in \A ,\, |x-\Re(L)| < |\Re(\epsilon)| \textrm{ and } \min_{n
% \in \Z} |y-\Im(L)
% + 2\pi n| \le |\Im(\epsilon)| & \textrm{ , if } n=3
%\end{array}\right\}.
%\end{displaymath}

We use the notation $f(t) \sim g(t)$ to mean that $\lim_{t \to \infty}
\frac{\Re(f(t))}{\Re(g(t))}=1$ and $\lim_{t \to \infty} \frac{\Im(f(t))}{\Im(g(t))}=1$.

%We think of a frame
%as an ordered orthonormal basis $f=(f_1,...,f_n)$. Throughout this
%paper, all frames will be positively oriented unless stated
%otherwise.

%define $\sim$.

%Define double cross, width, complex length of a geodesic. $G_L$.

%common perpendicular will refer to either the geodesic of the
%geodesic segment. It will often be oriented implicity by saying the
%common perpendicular from ...  to ....

%primitive and non-primitive geodesics.

%Check: when $n=2$, $\Omega_{2-1} = 2 \pi^{1}/\Gamma(2) = 2\pi$.

%When $n=3$, $\Omega_{3-1} = 2\pi^{3/2}/\Gamma(3/2)$. But
%$\Gamma(3/2)=(1/2)\Gamma(1/2) = (1/2)\sqrt{\pi}$. Thus,
%$\Omega_{3-1} = 4\pi$. This is correct.

%Let $\A$ denote either $\R$ or $\C/<2\pi i \Z>$ depending on whether $n=2$ or
%$3$.

%Given $L \in \A$ and $\epsilon \in \C$, let $B(L,\epsilon)$ be the set of all $a + ib \in
%\A$ such that $|a-\Re(L)| < \Re(\epsilon)$ and $\min_{n \in \Z} |b-\Im(L) + 2\pi n|<
%\Im(\epsilon)$.

\section{Counting Perpendiculars}

%Roughly speaking, given two geodesic segments in $M$ we will obtain
%an asymptotic formula (as $L \to \infty$) for the number of common
%perpendiculars between them with translation lengths in $B(L,\epsilon)$.

Let $\sigma_1=\sigma_1(t), \sigma_2=\sigma_2(t)$ be two oriented
geodesic segments in $M$. We allow them to depend on the parameter $t$
but we often suppress $t$ from the notation.

% For $i=1,2$, let $\sigma^f_i$ be a framing of
%$\sigma_i$ by which we mean that for each point $p$ in
%$\sigma_i$ ($i=1,2$) there is a unique frame $f(p) \in \sigma^f_i$ based at $p$ such that
%for any two points $p,q\in \sigma_i$, parallel transport along
%$\sigma_i$ from $p$ to $q$ maps $f(p)$ to $f(q)$. We also require
%that the first vector in the frame $f(p)$ is tangent to $\sigma_i$ and
%oriented consistently with the given orientation on $\sigma_i$.

For $i=1,2$, $T_1^\perp(\sigma_i)$ denote the set of unit vectors
perpendicular to $\sigma_i$. $T_1^\perp(\sigma_i)$ is naturally identified
with $\sigma_i \times S^{n-2}$. It carries the product volume
form in which the $\sigma_i$ factor has total volume
$length(\sigma_i)$ and the $S^{n-2}$ has total volume $\Omega_{n-2}$,
the volume of the unit $n-2$ sphere.

Any perpendicular $\gamma$ from $\sigma_1$ to $\sigma_2$ determines vectors $v_i(\gamma)
\in T_1^\perp(\sigma_i)$ in
the obvious way such that both $v_1(\gamma)$ and $v_2(\gamma)$ point toward $\gamma$.

Let $\chi=\chi_{\sigma_1,\sigma_2}$ denote the measure on
$T_1^\perp(\sigma_1)\times T_1^\perp(\sigma_2) \times \mF$ given by setting
$\chi(E)$ equal to the number of perpendiculars $\gamma$ from $\sigma_1$ to $\sigma_2$
such that
$(v_1(\gamma),v_2(\gamma),len(\gamma)) \in E$ where $len(\gamma)$
denotes the length of $\gamma$ if $n=2$ and $len(\gamma)$ denotes the
width of the double cross $(\sigma_1,\sigma_2;\gamma)$ otherwise. Here
$\gamma$ is oriented from $\sigma_1$ to $\sigma_2$.

Let $\chi'$ be the measure on
$T_1^\perp(\sigma_1)\times T_1^\perp(\sigma_2) \times \mF$ with derivative
\begin{eqnarray*}
d\chi' &=&  \frac{e^{(n-1)\Re(L)} }{2^{n-1}(2\pi)^{n-2} vol(T_1 M)}
dvol_{T_1^\perp(\sigma_1)} \, dvol_{T_1^\perp(\sigma_2)} \, dL
\end{eqnarray*}
where the measure on $\mF$ is standard Lebesgue measure.

Given a set $E \subset T_1^\perp(\sigma_1)\times T_1^\perp(\sigma_2) \times \mF$ and
$\epsilon >0$ there are two different $\epsilon$-neighborhoods of the boundary of $E$
that we will consider. The first is
\begin{eqnarray*}
N_1(\epsilon, \partial E)&:=& \{ (v_1,v_2,L) \in T_1^\perp(\sigma_1)\times
T_1^\perp(\sigma_2) \times \mF \, ; \, \exists (v'_1,v'_2,L) \in \partial E \\
&& \textrm{ such that } d(v'_1,v_1) \le \epsilon e^{-\Re(L)/2}, d(v'_2,v_2) \le \epsilon
e^{-\Re(L/2)}\}.
\end{eqnarray*}
In the above we used $d(\cdot,\cdot)$ to denote the usual distance in both
$T_1^\perp(\sigma_1)$ and $T_1^\perp(\sigma_2)$. The second is
\begin{eqnarray*}
N_2(\epsilon, \partial E) := \{ (v_1,v_2,L) \in T_1^\perp(\sigma_1)\times
T_1^\perp(\sigma_2) \times \mF ; \, \exists (v_1,v_2,L') \in \partial E \, \textrm{ and }
d(L,L') \le \epsilon \}.
\end{eqnarray*}

%Let $dg$ be the path metric on $T_1^\perp(\sigma_1)\times T_1^\perp(\sigma_2) \times
%\mF$ determined by
%\begin{eqnarray*}
%dg^2 = e^{\Re(L)/2}dv_1^2 + e^{\Re(L)/2}dv_2^2 + dL^2
%\end{eqnarray*}
%where $v_i \in T_1^\perp(\sigma_i)$, $L \in \mF$.

We say that a path $\{E_t\} \subset T_1^\perp(\sigma_1)\times T_1^\perp(\sigma_2) \times
\mF$ tends to infinity if
\begin{eqnarray*}
\inf \{\Re(l)\, | \, \exists (v_1,v_2,l) \in E_t\} \to \infty.
\end{eqnarray*}

\begin{thm}\label{thm:perp}
If $\{E_t\} \subset T_1^\perp(\sigma_1(t))\times T_1^\perp(\sigma_2(t)) \times \mF$ is a
path of compact sets tending to infinity and
\begin{eqnarray*}
\lim_{t \to \infty} \frac{\chi'(N_1(1, \partial E_t))}{\chi'(E_t)} =0 &\textrm{ and }&
\lim_{\epsilon \to 0} \limsup_{t \to \infty} \frac{\chi'(N_2(\epsilon, \partial
E_t))}{\chi'(E_t)} =0
\end{eqnarray*}
then $\chi(E_t) \sim \chi'(E_t)$.
\end{thm}

\begin{lem}\label{lem:reduction}
To prove theorem \ref{thm:perp} it suffices to prove the theorem for 'rectangular'
sequences $\{E_t\}$. These are sequences of the form
\begin{eqnarray*}
E_t=\sF_1 \times \sF_2 \times B(L,\epsilon) \subset T_1^\perp(\sigma_1)\times
T_1^\perp(\sigma_2) \times \mF
\end{eqnarray*}
where $\sF_1, \sF_2, L$ are functions of $t$, $\Re(\epsilon) >0$ if $n=2$, $\pi \ge
\Im(\epsilon)>0$ if $n=3$, $\sF_i=\sF'_i \times \sigma_i$, $\sF'_i$ is an arc of
$S^{n-2}$ with $length(\sF'_i)e^{L/2} \to \infty$ and $length(\sigma_i)e^{L/2} \to
\infty$.
\end{lem}

\begin{proof}
Let $\{E_t\}$ be a sequence satisfying the hypotheses of theorem \ref{thm:perp}. Let
$\pi: T_1^\perp(\sigma_1)\times
T_1^\perp(\sigma_2) \times \mF \to \mF$ be the projection map. Let $\epsilon >0$ if $n=2$
and $\Re(\epsilon), Im(\epsilon)>0$ if $n=3$. By partitioning $E_t$ if necessary and
reparametrizing we may assume that $\pi(E_t) \subset B(L,\epsilon)$. Because
\begin{eqnarray*}
\lim_{t \to \infty} \frac{\chi'(N_1(1, \partial E_t))}{\chi'(E_t)} =0
\end{eqnarray*}
there is a slowly increasing function $f(t)$ such that $f(t) \to \infty$ and
\begin{eqnarray*}
\lim_{t \to \infty} \frac{\chi'(N_1( f(t), \partial E_t))}{\chi'(E_t)} =0.
\end{eqnarray*}
Subdivide $T_1^\perp(\sigma_1)\times
T_1^\perp(\sigma_2)$ into rectangular solids with edge length $f(t)e^{-\Re(L)/2}$.
Subdivide $\mF$ into rectangles of edge length $\epsilon' >0$. Taking the product, we
obtain a subdivision of $T_1^\perp(\sigma_1)\times
T_1^\perp(\sigma_2) \times \mF$. Let $E'_t$ be the union of all those rectangles in the
subdivision which nontrivially intersect $E_t$. Because of the hypotheses on $E_t$,
\begin{eqnarray*}
\lim_{\epsilon' \to 0} \lim_{t \to \infty} \frac{\chi'(E'_t)}{\chi'(E_t)} = 1.
\end{eqnarray*}
So the conclusion of theorem \ref{thm:perp} is true for $\{E_t\}$ iff it is true for the
sequence $\{E'_t\}$. But each $E'_t$ is partioned into rectangles satisfying the
hypotheses of this lemma.
\end{proof}

For simplicity of exposition we will only prove the case
\begin{eqnarray*}
E_t = T_1^\perp(\sigma_1)\times T_1^\perp(\sigma_2) \times
B(L,\epsilon)
\end{eqnarray*}
the general case being similar. By choosing $\epsilon$ smaller if necessary we may assume
that $\Re(\epsilon)$ is less than the minimum injectivity radius of $M$.

By a vector pair $b=(x,y)$ we will mean that $x,y \in T_1(\H^n)$ (or in $T_1(M)$) and
that they share the same basepoint. We define an action of $Isom^+(\H^n)$ on the space of
vector pairs by setting $b \psi = (x',y')$ if there exists frames $f, f'$ with
$f=(x,y,...)$, $f'=(x',y',...)$ and $f \psi = f'$. The action of $Isom^+(\H^n)$ on the
frame bundle is discussed in section \ref{sec:note}. Define $-(x,y)=(-x,y)$.

Here is a proof sketch. We consider approximate perpendiculars obtained in the following
way. For
$i=1,2$, let $b_i=(x_i,y_i)$ be a vector pair where $x_i \in T_1^\perp(\sigma_i)$ and
$y_i$ has the same direction as $\sigma_i$. Consider $b_1G_{L/2}$ and $b_2G_{L/2}$. If
these two vector pairs
are close then there is a path $\gamma'$ formed by concatenating the
two paths traced out by $\{x_i G_s | 0 \le s \le \Re(L/2)\}$ ($i=1,
2$) with a short segment from the basepoint of $x_1G_{L/2}$ to that of
$x_2G_{L/2}$. Depending on how close $b_1G_{L/2}$ and $b_2G_{L/2}$ are, $\gamma'$ will be
close to a perpendicular $\gamma$ from $\sigma_1$ to $\sigma_2$. From this procedure we
obtain a function $P$ from the set of pairs $(b_1, b_2)$ of the above form to the set of
perpendiculars from $\sigma_1$ to $\sigma_2$. Asymptotically, the
length of $\gamma$ minus $L$ depends only on the position
of $b_1G_{L/2}$ relative to $b_2G_{L/2}$.

\begin{figure}[htb]
\begin{center}
\ \psfig{file=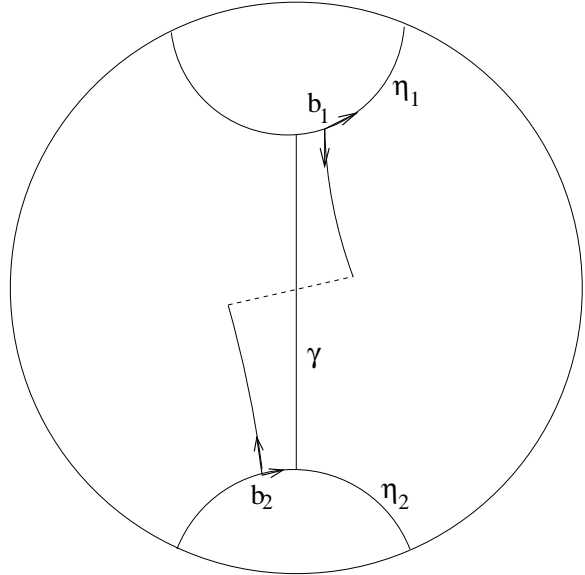,height=3in,width=3in}
\caption{The geodesics $\eta_1, \eta_2$ and the perpendicular $\gamma$ between them. The
'approximate perpendicular' determined by $b_1$ and $b_2$ is also shown.}
\label{fig:sketch}
\end{center}
\end{figure}

If we choose the pair $(b_1,b_2)$ uniformly at random, the distribution of
$(b_1G_{L/2},b_2G_{L/2})$ is asymptotically uniform for large $L$. To
obtain the asymptotics for $\chi(E_t)$ it suffices to estimate the volume of
$P^{-1}(\gamma)$ for each
perpendicular $\gamma$. This will be done by relating the calculations
to known results regarding the lattice point problem.

%The proof requires considering piecewise geodesic paths from $\sigma_1$ to $\sigma_2$
%that are close to being common perpendiculars. Later, we will relate the volume of the
%space of such paths to the number of common perpendiculars. It will then be important to
%understand the geometric relationship between such a path and the common perpendicular
%that it is close to. %

Let $\eta_1$ and $\eta_2$ be two disjoint oriented geodesics in $\H^n$. These will play
the role of lifts of $\sigma_1, \sigma_2$ to $\H^n$ in the sequel. Let $\gamma$ be the
perpendicular from $\eta_1$ to $\eta_2$. For $i=1,2$, let
$b_i$ be the vector pair $(x_i,y_i)$ where $x_i \in T_1^\perp(\eta_i)$ and $y_i$ points
in the direction of $\eta_i$.

%The piecewise geodesic path approximating $\gamma$ comes from connecting the segments
%traced out by $\{x_i G_t | 0 \le t \le \Re(L)\}$ (for $i=1,2$) with a short segment from
%the basepoint of $x_1G_{\Re(L)}$ to the basepoint of $x_2G_{\Re(L)}$.

Fix $L \in \mF$ and let $\phi \in Isom^+(\H^n)$ be such that
\begin{eqnarray}\label{eqn:b}
b_1 G_{L/2} \phi = -(b_2 G_{L/2}).
\end{eqnarray}
The pair $(L/2,\phi)$ determines the pair $(b_1,b_2)$ up to the
diagonal left action of $Isom^+(\H^n)$. Therefore, $(L/2,\phi)$ determines the oriented
geodesics $\eta_1,\eta_2$ and $\gamma$ up to rigid motion. By definition, $\eta_i$ is the
geodesic tangent to $y_i$ and $\gamma$ is the perpendicular from $\eta_1$ to $\eta_2$. It
now makes sense to define
$W_{L/2}(\phi)=width(\eta_1,\eta_2;\gamma)-L$. In other
words $W_{L/2}(\phi)$ is the error in the ``length'' of the perpendicular corresponding to
$(b_1,b_2)$
when equation \ref{eqn:b} is satisfied.

If $v$ is a unit vector on $\H^n$, let $H(v)$ be the unique horoball for which $v$ is an
outer unit normal.  Let $\Phi \subset Isom^+(\H^n)$ be the set of all isometries $\phi$
such
that for each unit vector $v$ on $\H^n$, $H(v)$ and $H(-v\phi)$ intersect trivially.

\begin{lem}\label{lem:w}
There exists a continuous function $W_\infty: \Phi \to \mF$ such that
\begin{eqnarray*}
\lim_{L\to \infty} W_{L/2}(\phi) &=& W_\infty(\phi).
\end{eqnarray*}
\end{lem}

\begin{proof}
Fix $\phi \in \Phi$. For $i=1,2$, let $b'_i=(x'_i,y'_i)$ be a vector pair on $\H^n$ so
that
\begin{eqnarray*}
b'_1 \phi = -b'_2.
\end{eqnarray*}
For $L \in \mF$, let $b_i(L/2)=b'_i G_{-L/2}$. Equation \ref{eqn:b} is
satisfied. So let $\eta_1=\eta_1(L/2), \eta_2=\eta_2(L/2), \gamma=\gamma(L/2)$ be as in
the paragraphs preceding this lemma. That is, if $b_i=(x_i,y_i)$ then $\eta_i$ is tangent
to $y_i$ and $\gamma$ is perpendicular to $\eta_1$ and $\eta_2$. Let $C_i(L/2)$ be the
radius $\Re(L/2)$ cylinders
surrounding $\eta_i(L/2)$. The basepoint of $b'_i$ is in $C_i(L/2)$. Let $\alpha_i(L/2)$
be a geodesic perpendicular to $\gamma$
which passes through the intersection point $\gamma \cap C_i(L/2)$ such that
$width(\eta_i,\alpha_i;\gamma)=(-1)^{i-1} L/2$. Then
\begin{eqnarray*}
width(\alpha_1,\alpha_2;\gamma)&=&width(\eta_1,\eta_2;\gamma) -
width(\eta_1,\alpha_1;\gamma)-width(\alpha_2,\eta_2; \gamma)\\
&=& width(\eta_1,\eta_2;\gamma)-L = W_{L/2}(\phi).
\end{eqnarray*}

%\begin{figure}[htb]
%\begin{center}
%\ \psfig{file=cylinders2.eps,height=2.5in,width=2.5in} \caption{The
%cylinders $C_1(L),C_2(L)$ surround $\eta_1(L),\eta_2(L)$. The
%first two components of the frames $\hf^1,\hf^2, h^1(L), h^2(L)$ are
%shown.} \label{fig:cylinders}
%\end{center}
%\end{figure}

We finish the lemma by taking the limit of both sides as $L\to \infty$. The unit vector
parallel to $\alpha_i$ at its intersection
point with $\gamma$ is parallel to $y'_i$ with respect to the Euclidean
structure on the cylinder $C_i$. As $L$ tends to infinity, $C_i(L/2)$ tends to the
horosphere
$H_i$ for which $x'_i$ is an outer unit normal. Since $\phi \in
\Phi$, $H_1$ and $H_2$ intersect trivially. Since $\gamma(L/2)$ is
at right angles to $C_1(L/2)$ and $C_2(L/2)$, $\gamma(L/2)$ converges to the unique
geodesic $\gamma(\infty)$ intersecting both $H_1$
and $H_2$ at right angles. $\alpha_i(L/2)$ converges to the unique geodesic perpendicular
to $\gamma(\infty)$ that passes through the intersection point $\gamma(\infty) \cap
H_i(L/2)$ and is parallel to $y'_i$ in $H_i$.

Thus $W_{L/2}(\phi)$ converges to $width(\alpha_1(\infty), \alpha_2(\infty);
\gamma(\infty))=:W_\infty(\phi)$. Since this double cross depends only on $b'_1$ and
$b'_2$, it is
independent of $L$. Continuous dependence on $\phi$ is clear.
\end{proof}

Let $\Phi_{L/2}=\Phi_{L/2}(\epsilon)$ be the set of isometries $\phi \in \Phi$ such that
$W_{L/2}(\phi) \in
B(0,\epsilon)$ and the distance between $v$ and $v\phi$ is at most $\Re(\epsilon)$ for any
vector $v \in T_1(\H^n)$. Let $\Phi_\infty=\Phi_\infty(\epsilon)$ be the set of all $\phi
\in \Phi$ such that
$W_\infty(\phi) \in B(0,\epsilon)$ and $d(v,v\phi)\le \Re(\epsilon)$ for any vector $v$
on $\H^n$. By
the above lemma, $\Phi_{L/2}$ converges to $\Phi_\infty$ in the Hausdorff topology.

%satisfying:
%\begin{itemize}
%\item $W(L,\phi) \in B(0,\epsilon)$,
%\item the distance between the basepoint of $f^0$ and and the basepoint of $\phi f^0$ is
%at most $|\epsilon|$.
%\end{itemize}

%Given two fixed geodesics $\eta_1, \eta_2$, we will need to understand the volume of the
%space of pairs $(b_1,b_2)$ satisfying $b_1G_{L/2}\phi = - (b_2G_{L/2})$ where $\phi \in
%\Phi_{L/2}$. This will later allow us
%to obtain the number of perpendiculars between $\sigma_1, \sigma_2$ from the volume of
%the space of the approximate perpendiculars.

For a fixed $\phi \in \Phi$ and $L\in \mF$, suppose $(b_1,b_2)$ satisfies equation
\ref{eqn:b}. Let $\eta_1, \eta_2, \gamma$ be as above (i.e. $\eta_i$ is the geodesic
tangent to $y_i$ and $\gamma$ is perpendicular to $\eta_1$ and $\eta_2$). Let
$A_{L/2}(\phi)$ be the set of
all $(\bb_1,\bb_2)$ where for $i=1,2$, $\bb_i=(\bx_i,
\by_i)$ is a vector pair and
\begin{itemize}
\item $\bx_i \in T_1^\perp(\eta_i)$,
\item $\by_i$ points in the direction of $\eta_i$,
\item there exists $\psi \in \Phi_{L/2}$ with
\begin{eqnarray}\label{eqn:z}
\bb_1 G_{L/2} \psi = -(\bb_2 G_{L/2}).
\end{eqnarray}
\end{itemize}
We can identify
$A_{L/2}(\phi)$ with a subset of $T_1^\perp(\eta_1)\times
T_1^\perp(\sigma_2)$ by the map $(\bb_1,\bb_2)\to(\bx_1,\bx_2)$. This
map is 1-1 since $\by_i$ is uniquely determined by $\bx_i$. By definition, the volume of
$A_{L/2}(\phi)$ is its volume in $T_1^\perp(\eta_1)\times
T_1^\perp(\eta_2)$ with respect to product measure.

%We will let $vol(\cdot)$ denote the usual volume form on
%$T_1^\perp(\eta_i)$ or $T_1^\perp(\eta_1)\times T_1^\perp(\eta_2)$ or
%on the unit tangent bundles $T_1(C_i)$ or on $T_1(C_1) \times T_1(C_2)$
%depending on the context. Applied to a product space, $vol(\cdot)$
%will mean the usual product measure. Applied to $T_1^\perp(\eta_i)$,
%$vol(\cdot)$ equals $\nu_i$ times $\Omega_{n-2}$, the volume of the
%unit $n-2$ sphere.

\begin{lem}\label{lem:cylinder volume}
There exists a continuous function $V$ on $\Phi_\infty$ such that
\begin{eqnarray*}
V(\phi)&=&\lim_{L \to \infty} vol(A_{L/2}(\phi)) \cosh(\Re(L)/2)^{2n-2}.
\end{eqnarray*}
Also the function $\phi \to vol(A_{L/2}(\phi)) \cosh(\Re(L)/2)^{2n-2}$ converges to $V$
in $L^1(\Phi_\infty)$ as $L\to \infty$.
\end{lem}
%{\bf Remark:} The $n$-dimensional analogue of the above formula is $vol(A_L(\phi)) =
%v_L(\phi) \cosh(L)^{2-2n}$.

\begin{proof}
Fix $\phi \in \Phi$. Let $b'_i, b_i, \eta_i, \gamma, C_i, \alpha_i, H_i$ be as in the
previous lemma. If for some $i=1,2$, two unit vectors $v_0,v_1 \in T_1^\perp(\eta_i)$
share the same basepoint, then the
distance between the basepoints of $v_0G_t$ and $v_1G_t$ in $C_i(t)$ is $\sinh(t)$ times
the angle between them. If two vectors $v_0, v_1 \in T_1^\perp(\eta_i)$ are parallel
along $\eta_i$ then the distance between the basepoints of $v_0G_t$ and $v_1G_t$ in
$C_i(t)$ is $\cosh(t)$ times the distance between the basepoints of $v_0$ and $v_1$.
Therefore, if $X \subset T_1^\perp(\eta_i)$ then
\begin{eqnarray*}
vol(XG_{\Re(L/2)}) &=&  \cosh(\Re(L/2))^{n-1} vol(X) + O(vol(X)e^{-(n-1)\Re(L/2)}).
\end{eqnarray*}
Thus
\begin{eqnarray*}
vol(A_{L/2}(\phi) G_{\Re(L/2)}) &=&  \cosh(\Re(L/2))^{2n-2} vol(A_{L/2}(\phi)) +
O(vol^2(A_{L/2}(\phi))).
\end{eqnarray*}
Here $A_{L/2}(\phi)G_{\Re(L/2)}$ refers to the diagonal action of $G_t$ on
$T_1^\perp(\eta_1)\times T_1^\perp(\eta_2)$.

As $L$ tends to infinity, the cylinders $C_i(L/2)$ tend to the horospheres $H_i$. Let
$T_1^\perp(H_i)$ be the set of outward unit normals to $H_i$. Let $A_\infty(\phi)$ be the
set of pairs $(a_1,a_2) \in T_1^\perp(H_1) times T_1^\perp(H_2)$ such that there exists
$\psi \in \Phi_\infty$ with $a_1 \psi = -a_2$.

By
definition, a pair $(a_1,a_2) \in A_{L/2}(\phi)G_{\Re(L/2)}$ if and only if there exists
$\psi \in \Phi_{L/2}$ such that $a_1 \psi = -a_2$. So $A_{L/2}(\phi)G_{\Re(L/2)}$
converges in the Hausdorff topology to $A_\infty(\phi)$. Because $A_{L/2}(\phi)$ is
defined by a finite set of real analytic inequalities, the volume of
$A_{L/2}(\phi)G_{\Re(L/2)}$ converges to the volume of $A_\infty(\phi)$. Set
$V(\phi):=vol(A_\infty(\phi))$. It is clear that $V$ is continous. $L^1$-convergence
follows from pointwise convergence and the boundedness of the two functions.
\end{proof}

For $i=1,2$, let $T_2^\perp(\sigma_i)$ denote the space of vector pairs
$b=(x,y)$ such that $x \in T_1^\perp(\sigma_i)$ and $y$ has the same direction as
$\sigma_i$. $T_2^\perp(\sigma_i)$ is naturally
identified with $T_1^\perp(\sigma_i)$ by the map $b \to x$.

Let $T_2(M)$ denote the space of all vector pairs
$(x,y)$ where $x, y$ are unit vectors in $M$ which share the same
basepoint and are orthogonal to each other. $T_2(M)$ splits locally as
a product $M \times S^{n-1} \times S^{n-2}$. It
has a natural volume form induced from its Riemannian structure with
total mass $vol(M)\Omega_{n-1}\Omega_{n-2}$. The standard probability measure on $T_2(M)$
is this volume form normalized to have total mass $1$.

\begin{lem}\label{lem:equidistribution}
For $i=1,2$ let $\mu_i=\mu_{i,t}$ be the probability measure on
$T_2^\perp(\sigma_i(t))$ obtained from normalizing the volume form on
$T_1^\perp(\sigma_i(t))$. If $\cosh(\Re( L(t)) /2)length(\sigma_i(t))\to \infty $ and
$L(t) \to \infty$ as $t \to
\infty$ then the pushforward measure $G_{\Re(L/2)*}(\mu_i)$ converges
in the weak* topology to the standard probability measure on $T_2(M)$.
\end{lem}
\begin{proof}
Let $\lambda_i$ be a subsequential weak* limit of
$\{(G_{L/2})_*(\mu_i)\}_t$. Let $\tsigma_i$ be a lift of $\sigma_i$ to
the universal cover $\H^n$. Let $C$ be the radius $\Re(L/2)$ cylinder with axis
$\tsigma_i$. Then $T_2^\perp(\tsigma_i) G_{\Re(L/2)}$ is the set of
vector pairs $b=(x,y)$ with $x$ an outward unit normal to $C$ and $y$
parallel to the axis of $C$. The length of this cylinder is $\cosh(L/2)
length(\sigma_i)$ which tends to infinity with $L$ by assumption. So
for any radius $R$, the probability that a point chosen uniformly at
random on $C$ is within distance $R$ of the boundary of the cylinder
tends to zero as $L$ tends to infinity. As noted before, the
geometry of the cylinder tends towards that of a horosphere with
increasing $L$. It follows that $\lambda_i$ is invariant under a
rank $n-1$ parabolic subgroup of $Isom^+(\H^n)$. Ratner's theorems \cite{Ratner} on
unipotent
flows now
imply the lemma.
\end{proof}

\begin{lem}\label{lem:rough count}
Suppose that for $i=1,2$, $length(\sigma_i(t))\cosh(L(t)/2) \to \infty$ as
$t \to \infty$. Then
\begin{eqnarray*}
\chi\left(T_1^\perp(\sigma_1) \times T_1^\perp(\sigma_2) \times B(L,\epsilon)\right) \sim
c_n(\epsilon) e^{(n-1)\Re(L)}
\frac{length(\sigma_1)length(\sigma_2)}{vol(FM)}
\end{eqnarray*}
where
\begin{eqnarray*}
c_n(\epsilon) = 2^{2-2n} \Omega_{n-2}^2 \int_{\Phi_\infty} \frac{1}{V(\phi)} d\phi.
\end{eqnarray*}
The above integral formula is with respect to Haar measure on $Isom^+(\H^n)$.
\end{lem}

%{\bf Remark }: Since the above formula depends on $\sigma_1$ and $\sigma_2$ only though
%the product of their lengths, it follows that if a perpendicular
%$\gamma \in \sP$ is chosen uniformly at random, then its endpoints
%$(e_1,e_2)$ in $(\sigma_1,\sigma_2)$ are (asymptotically)
%distributed uniformly at random.

\begin{proof}
Let $\Delta_{L/2}$ be the set of all pairs $(b_1,b_2)$ in which $b_i =(x_i,y_i)\in
T_2^\perp(\sigma_i)$ and there exists $\phi \in \Phi_{L/2}$ such that
\begin{eqnarray}\label{eqn:b2}
b_1 G_{L/2} \phi = -(b_2 G_{L/2}).
\end{eqnarray}

Associated to every such pair $(b_1, b_2)$ is a geodesic segment $\gamma$ perpendicular
to $\sigma_1$ and $\sigma_2$ at its endpoints defined as follows.

The distance from $x_1G_{L/2}$ to $x_2G_{L/2}$ is at most $\Re(\epsilon)$ since $\phi \in
\Phi_{L/2}(\epsilon)$. Let $\gamma'$ be the path formed from concatening the
segments formed from the trace of $\{x_i G_t : 0 \le t \le L/2\}$
(for $i=1,2$) with a segment of length at most $\Re(\epsilon)$ joining
the basepoint of $x_1G_{L/2}$ to that of $x_2G_{L/2}$. The middle segment is unique since
by hypothesis the minimum injectivity radius of $M$ is at least $\Re(\epsilon)$. Let
$\tgamma'$ be a lift of $\gamma'$ to the universal cover and let
$\tsigma_1, \tsigma_2$ be the lifts of $\sigma_1, \sigma_2$ that
pass through the endpoints of $\tgamma'$. Let $\tgamma$ be the common
perpendicular to $\tsigma_1$ and $\tsigma_2$ and $\gamma$ its
projection to $M$.

The definition of $\Phi_{L/2}(\epsilon)$ implies $len(\gamma)\in B(L,\epsilon)$. Define
$P(b_1,b_2)=\gamma$. Conversely, if $\gamma$ is a segment perpendicular to
$\sigma_1$ and $\sigma_2$ at its
endpoints and $len(\gamma) \in B(L,\epsilon)$ then there exists $(b_1,b_2) \in
\Delta_{L/2}$ such that $P(b_1,b_2)=\gamma$. To see this, let $b_i=(x_i,y_i) \in
T_2^\perp(\sigma_i)$ where
$x_i=v_i(\gamma)$ is as defined in the paragraphs before theorem \ref{thm:perp}. Then
$\gamma'=\gamma$ and $P(b_1,b_2)=\gamma$.

Let $(b_1,b_2) \in \Delta_{L/2}$ be chosen uniformly at random. Let
$\phi \in \Phi_{L/2}$ be such that equation \ref{eqn:b2} holds. Let $V(b_1,b_2)$ denote
the volume
of
\begin{eqnarray*}
\{ (b'_1,b'_2) \in \Delta_{L/2} : P(b_1,b_2)=P(b'_1,b'_2)\}.
\end{eqnarray*}
So $V(b_1,b_2) = vol(A_{L/2}(\phi))$ where $A_{L/2}(\cdot)$ is defined
as before lemma \ref{lem:cylinder volume} with $\tsigma_i$ replacing $\eta_i$. Note
\begin{eqnarray*}
\chi\left(T_1^\perp(\sigma_1)\times T_1^\perp(\sigma_2) \times B(L,\epsilon)\right) &=&
vol(\Delta_{L/2}) \mE[1/V(b_1,b_2)]
\end{eqnarray*}
where $\mE$ denotes expected value. We claim
\begin{eqnarray*}
vol(\Delta_{L/2}) \sim \Omega^2_{n-2}
length(\sigma_1)length(\sigma_2)\frac{vol(\Phi_\infty)}{
vol(FM)}.
\end{eqnarray*}
To see this note that the total mass of $T_2^\perp(\sigma_1)\times T_2^\perp(\sigma_2)$
is $\Omega^2_{n-2}length(\sigma_1)length(\sigma_2)$. So

 $$vol(\Delta_{L/2})/ \Omega^2_{n-2}
length(\sigma_1)length(\sigma_2)$$
is the probability that a pair
$(b_1,b_2)$ chosen uniformly at random in $T_2^\perp(\sigma_1)\times
T_2^\perp(\sigma_2)$ is in $\Delta_{L/2}$. By lemma \ref{lem:equidistribution} that is
asympotic to the probability that a uniformly random pair $(b'_1,b'_2) \in T_2 M \times
T_2 M$ satisfies $b'_1\psi = b'_2$ for some $\psi \in \Phi_{L/2}(\epsilon)$. The latter
equals $vol(\Phi_{L/2})/vol(FM)$ which is asymptotic to $vol(\Phi_\infty)/vol(FM)$.

%This choice is unique only up to
%precomposition with an element of $SO(n-2)$ since
%this does not change $b_1G_{L/2}$. Since $A_{L/2}(\phi)$ is invariant under $SO(n-2)$, it
%follows from lemma
%\ref{lem:cylinder volume} that $v$ is $SO(n-2)$ invariant too.

By lemma \ref{lem:cylinder volume}, $V(b_1,b_2)$ is asymptotic to
$V(\phi)\cosh^{2-2n}(\Re(L)/2)$ where $\phi \in Isom^+(\H^n)$ is any
element such that equation \ref{eqn:b2} holds. Lemma \ref{lem:equidistribution} implies
\begin{eqnarray*}
\mE[1/V(b_1,b_2)] &\sim & \frac{1}{vol(\Phi_\infty)}
\int_{\Phi_\infty} \frac{1}{V(\phi)\cosh^{2-2n}(\Re(L)/2)} d\phi.
\end{eqnarray*}
So,
\begin{eqnarray*}
\chi(T_1^\perp(\sigma_1)\times T_1^\perp(\sigma_2) \times
B(L,\epsilon)) &=&vol(\Delta_{L/2}) \mE[1/V(b_1,b_2)]\\
%&\sim& \Omega_{n-2}^2 length(\sigma_1)length(\sigma_2)\frac{vol(\Phi_\infty)}{vol(FM)}
%\\
%&& \times \frac{1}{vol(\Phi_\infty)}
%\int_{\Phi_\infty}
%\frac{1}{v(\phi)\cosh^{2-2n}(\Re(L)/2)} d\phi\\
&\sim& c_n(\epsilon) e^{(n-1)\Re(L)}
\frac{length(\sigma_1)length(\sigma_2)}{vol(FM)}
\end{eqnarray*}
where
\begin{eqnarray*}
c_n(\epsilon) = 2^{2-2n} \Omega_{n-2}^2 \int_{\Phi_\infty} \frac{1}{V(\phi)} d\phi.
\end{eqnarray*}

\end{proof}

The next step in the proof of theorem \ref{thm:perp} is to compute
$c_n(\epsilon)$. The previous lemmata can all be
modified slightly by replacing the segments $\eta_i$ and $\sigma_i$
with points. The resulting lemmas apply to the lattice point problem. The value of
$c_n(\epsilon)$ will then follow from classical results.

For $L, \epsilon_1$ positive real numbers and $p_1,p_2 \in M$, let
$N(p_1,p_2,L,\epsilon_1)$ be the set of
geodesic segments from $p_1$ to $p_2$ with length in the interval
$(L-\epsilon_1,L+\epsilon_1)$. The analogue of lemma \ref{lem:rough count} is:

\begin{lem}\label{lem':rough count}
\begin{eqnarray*}
|N(p_1,p_2,L,\epsilon_1)| \sim c'_n(\epsilon_1 + i\pi) e^{(n-1)L}{vol(M)}
\end{eqnarray*}
where
\begin{eqnarray*}
c'_n(\epsilon) &=& \frac{\Omega_{n-1}^2 c_n(\epsilon)}{\Omega_{n-2}^2}.
\end{eqnarray*}
\end{lem}

Before proving the above, we need analogues for lemmas \ref{lem:w},
\ref{lem:cylinder volume} and \ref{lem:equidistribution}. Let $v_1, v_2$ be unit vectors
on $\H^n$. Suppose $\phi \in Isom^+(\H^n)$ satisfies
\begin{eqnarray}\label{eqn:v}
v_1 G_{L/2} \phi = -(v_2 G_{L/2}).
\end{eqnarray}
The pair $(v_1,v_2)$ is determined up to a rigid motion by $L, \phi$
and the above equation. So it makes sense to define $W'_{L/2}(\phi)$ to be the
distance between the basepoints of $v_1$ and $v_2$ minus $L$.

\begin{lem}\label{lem':w}
If $\phi \in \Phi$ then
\begin{eqnarray*}
\lim_{L\to \infty} W'_{L/2}(\phi) &=& \Re(W_\infty(\phi)).
\end{eqnarray*}
where $W_\infty$ is as defined in lemma \ref{lem:w}.
\end{lem}

\begin{proof}
The proof is similar to the proof of lemma \ref{lem:w}. The cylinders
$C_i(L/2)$ are replaced with spheres $S_i(L/2)$ of radius $L/2$ centered at the
basepoint of $v_i$.
\end{proof}

Let $\Phi'_{L/2}=\Phi'_{L/2}(\epsilon)$ be the set of isometries $\phi \in \Phi$
satisfying:
\begin{itemize}
\item $W'_{L/2}(\phi) \in (-\epsilon_1, \epsilon_1)$,
\item for any vector $v \in T_1(\H^n)$, the distance between $v$ and
$v\phi$ is at most $\epsilon$.
\end{itemize}
As before, $\Phi_\infty$ is defined as above except that the first condition is
replaced by $W_\infty(\phi) \in (-\epsilon_1,\epsilon_1)$. Lemma \ref{lem':w}
implies that $\Phi'_{L/2}$ converges to $\Phi_\infty$ in the Hausdorff
topology.

For a fixed $\phi \in \Phi$ and $L\in \mF$, suppose $(v_1,v_2)$
satisfies equation \ref{eqn:v}. For $i=1,2$, let $q_i$ be the basepoint
of $v_i$. Let $A'_{L/2}(\phi)$ be the set of all $(z_1,z_2) \in T_1(q_1) \times T_1(q_2)$
such that there exists $\psi \in \Phi'_{L/2}$ with
\begin{eqnarray}\label{eqn':z}
z_1 G_{L/2} \psi = -(z_2 G_{L/2}).
\end{eqnarray}

\begin{lem}\label{lem':cylinder volume}
With $V:\Phi_\infty \to \R$ defined as in lemma \ref{lem:cylinder volume},
\begin{eqnarray*}
V(\phi)&=&\lim_{L \to \infty} vol(A'_{L/2}(\phi)) \cosh(L/2)^{2n-2}
\end{eqnarray*}
for all $\phi \in \Phi_\infty$. Convergence also holds in $L^1(\Phi_\infty)$.
\end{lem}

\begin{proof}
The proof is similar to the proof of \ref{lem:cylinder volume}.
\end{proof}

\begin{lem}\label{lem':equidistribution}
For $i=1,2$ let $\mu'_i$ be the obvious probability measure on
$T_1(q_i)$. The pushforward measure $\{G_{L/2 *}(\mu_i)\}_L$ converges
in the weak* topology to the standard probability measure on $T_1 M$ obtained by
normalizing its volume form.
\end{lem}

\begin{proof}
This result is classical. See [Bekka and Mayer, 2000] for example.
\end{proof}

%The proof of lemma \ref{lem':rough count} now follows almost line for
%line like the proof of lemma \ref{lem:rough count}. We leave the
%details to the reader.

\begin{proof}(of lemma \ref{lem':rough count})
The proof is similar to the proof of lemma \ref{lem:rough count}. Let $\Delta'_{L/2}$ be
the set of
all pairs $(v_1,v_2)\in T_1(p_1)\times T_1(p_2)$ such that there exists $\phi \in
\Phi'_{L/2}$ such that
\begin{eqnarray}\label{eqn':b2}
v_1 G_{L/2} \phi = -(v_2 G_{L/2}).
\end{eqnarray}
For $(v_1,v_2) \in \Delta'_{L/2}$ let $P(v_1,v_2)$ denote the
corresponding segment
from $p_1$ to $p_2$. Let $V(v_1,v_2)$ denote the volume of the set of
all $(v'_1,v'_2) \in \Delta'_{L/2}$ with $P(v'_1,v'_2)=P(v_1,v_2)$. As
in lemma \ref{lem:rough count},
\begin{eqnarray*}
|N(p_1,p_2,L,\epsilon_1)| &=& vol(\Delta'_{L/2}) \mE[1/V(v_1,v_2)]
\end{eqnarray*}
where $\mE$ denotes expected value. Since $vol(T_1(p_1)\times T_1(p_2)) = \Omega_{n-1}^2$
it follows from
lemma \ref{lem':equidistribution} that
\begin{eqnarray*}
vol(\Delta'_{L/2}) \sim \frac{\Omega^2_{n-1} vol(\Phi_\infty)}{vol(FM)}.
\end{eqnarray*}
%Here we have made used of the normalizations $vol(M)=vol(T_1 M)$ and
%$length(\sigma_i)=vol(T_1^\perp(\sigma_i))$.

By lemma \ref{lem':cylinder volume}, $V(v_1,v_2)$ is asymptotic to
$V(\phi)\cosh^{2-2n}(\Re(L)/2)$ where $\phi \in Isom^+(\H^n)$ is any
element satisfying equation \ref{eqn':b2}. So
\begin{eqnarray*}
\mE[1/V(v_1,v_2)] &\sim & \frac{1}{vol(\Phi_\infty)}
\int_{\Phi_\infty} \frac{1}{V((\phi))\cosh^{2-2n}(\Re(L)/2)} d\phi.
\end{eqnarray*}
So
\begin{eqnarray*}
|N(p_1,p_2,L,\epsilon_1)| &=& vol(\Delta'_{L/2}) \mE[1/V(b_1,b_2)]\\
%&\sim& \frac{vol(\Phi_\infty)\Omega^2_{n-1}}{vol(FM)} \frac{1}{vol(\Phi_\infty)}
%\int_{\Phi_\infty}
%\frac{1}{v(\phi)\cosh^{2-2n}(\Re(L)/2)} d\phi\\
&\sim& c'_n(\epsilon) e^{(n-1)\Re(L)}
\frac{1}{vol(FM)}
\end{eqnarray*}
where
\begin{eqnarray*}
c'_n(\epsilon) &=& 2^{2-2n} \Omega_{n-1}^2 \int_{\Phi_\infty}
\frac{1}{V(\phi)} d\phi = \frac{\Omega_{n-1}^2 c_n(\epsilon)}{\Omega_{n-2}^2}.
\end{eqnarray*}
\end{proof}

\begin{lem}\label{lem:c}
If $n=2$, let $\epsilon>0$. Otherwise let $\epsilon = \epsilon_1 + i\epsilon_2 \in \A$
with $\epsilon_1 >0$ and $0 < \epsilon_2 \le \pi$. Then
\begin{displaymath}
c_n(\epsilon) = \left\{\begin{array}{ll}
\frac{\Omega^2_{n-2}\sinh(\epsilon)vol(FM)}{\Omega_{n-1}(n-1)2^{n-2}vol(M)} & \textrm{
if } n=2\\
\frac{2\epsilon_2}{2\pi}\frac{\Omega^2_{n-2}\sinh((n-1)\epsilon_1)vol(FM)}{\Omega_{n-1}(n-1)2^{n-2}vol(M)} &
\textrm{
if
}
n=3.
\end{array}\right\}
\end{displaymath}
\end{lem}
\begin{proof}
It is a classical result (Bekka and Mayer, 2000) that
\begin{eqnarray*}
|N(p_1,p_2,L,\epsilon_1)| &\sim& \frac{vol(B_{L+\epsilon_1}) -
vol(B_{L-\epsilon_1})}{vol(M)}
\end{eqnarray*}
when $L>0$. The volume of $B^n_r$, the radius $r$ ball in $\H^n$, is
\begin{eqnarray*}
vol(B^n_r) &=& \Omega_{n-1} \int_0^r \sinh^{n-1}(t) \, dt \sim  \Omega_{n-1}
\frac{\cosh^{n-1}(r)}{n-1}.
\end{eqnarray*}
So
\begin{eqnarray*}
|N(p_1,p_2,L,\epsilon_1)| &\sim& \frac{\Omega_{n-1}
\sinh((n-1)\epsilon_1)e^{(n-1)L}}{(n-1)2^{n-2}vol(M)}.
\end{eqnarray*}
By the previous lemma,
\begin{eqnarray*}
|N(p_1,p_2,L,\epsilon_1)| &\sim&  \frac{c'_n(\epsilon_1 + i\pi) e^{(n-1)\Re(L)}
}{vol(FM)}.
\end{eqnarray*}
Therefore,
\begin{eqnarray*}
c'_n(\epsilon_1 + i\pi)
&=&\frac{\Omega_{n-1}\sinh((n-1)\epsilon_1)vol(FM)}{(n-1)2^{n-2}vol(M)}.
\end{eqnarray*}
Now,
\begin{eqnarray*}
c_n(\epsilon) &=& \frac{\Omega_{n-2}^2}{\Omega_{n-1}^2}
c'_n(\epsilon) =
\frac{\Omega^2_{n-2}\sinh((n-1)\epsilon_1)vol(FM)}{\Omega_{n-1}(n-1)2^{n-2}vol(M)}.
\end{eqnarray*}
This proves the cases $n=2$ and $\epsilon_2=\pi$. Lemma \ref{lem:rough count} implies that
\begin{eqnarray*}
\chi\left( T_1^\perp(\sigma_1) \times T_1^\perp(\sigma_2) \times B(L,\epsilon)\right)
&\sim& \chi\left( T_1^\perp(\sigma_1) \times T_1^\perp(\sigma_2) \times B(L +
ix,\epsilon)\right)
\end{eqnarray*}
for any $x \in \R$. So if $\gamma=\gamma_L$ is a common
perpendicular between $\sigma_1$ and $\sigma_2$ chosen uniformly at random among all
perpendiculars with $|\Re(len(\gamma))-\Re(L)| < \epsilon$ then the distribution of
$\Im(len(\gamma))$ is asymptotically uniform
in $[-\pi, \pi]$. This implies equidistribution of holonomy which finishes the lemma.
\end{proof}

\begin{proof}(of theorem \ref{thm:perp})
Assume $\{E_t\}$ is rectangular in the sense of lemma \ref{lem:reduction}. In the case
when $\sF_1=T_1^\perp(\sigma_1), \sF_2=T_1^\perp(\sigma_2)$ and $n =3$ it
follows from  lemma \ref{lem:rough count} and lemma \ref{lem:c} above that
\begin{eqnarray*}
\chi\left( E_t\right) &\sim&
\frac{2|\Im(\epsilon)|}{2\pi}\frac{\Omega^2_{n-2}\sinh(\Re((n-1)\epsilon))vol(FM)}{\Omega_{n-1}(n-1)2^{n-2}vol(M)}e^{(n-1)\Re(L)}
\frac{length(\sigma_1)length(\sigma_2)}{vol(FM)}\\
&=& \frac{2|\Im(\epsilon)|}{2\pi}\frac{\sinh(\Re((n-1)\epsilon))}{(n-1)2^{n-2}vol(T_1 M)}
e^{(n-1)\Re(L)} vol(T_1^\perp(\sigma_1))vol(T_1^\perp(\sigma_2)).
\end{eqnarray*}
Here we used that $vol(T_1^\perp(\sigma_i)) = \Omega_{n-2}length(\sigma_i)$ and $vol(T_1
M)=\Omega_{n-1}vol(M)$. A routine calculation shows that this equals $\chi'(E_t)$. The
$n=2$ case is similar. The general case is also similar and is left to the reader.

% On the other hand,
%\begin{eqnarray*}
%\chi'(E_L) &=& \frac{1}{(2\pi)^{n-2} 2^{n-1} vol(T_1 M)}\int_{E_L} \,
%e^{(n-1)\Re(L)} dvol_{T_1^\perp(\sigma_1)} \, dvol_{T_1^\perp(\sigma_2)} \, dL\\
%&=& \frac{2|\Im(\epsilon)|}{2\pi}\frac{\sinh(\Re(\epsilon))}{(n-1)2^{n-2}vol(T_1 M)}
%e^{(n-1)\Re(L)} vol(T_1^\perp(\sigma_1))vol(T_1^\perp(\sigma_2)).
%\end{eqnarray*}

\end{proof}

\section{Pants and Free Subgroups: Definitions and Notation}\label{sec:pants}
We need more notation for pants. Consider the dimension 2 case first. Let $P$ be a pair
of pants with a hyperbolic metric and geodesic boundary components $\partial_1,
\partial_2, \partial_3$. Let $l_i$ be the length of $\partial_i$. For $i \ne j$ define
$l_{ij}(l_1,l_2,l_3)$ to be the length of the shortest path from $\partial_i$ to
$\partial_j$. We will give a formula for $l_{12}$ that is asymptotic as all three lengths
$l_1,l_2,l_3\to\infty$.

In the case of dimension 3, consider a pair of pants $P$ with oriented boundary
components $\partial_i$ ($i=1,2,3$). Let $\phi:\pi_1(P) \to Isom^+(\H^3)$ be a
homomorphism. Let * be a basepoint in the interior of $P$. Let $\partial'_i$ be a loop
based at * freely homotopic to $\partial_i$. By choosing $\partial'_i$ appropriately we
may assume that $[\partial'_1]$ and $[\partial'_2]$ freely generate $\pi_1(P,*)$ and
$[\partial'_3]=[\partial'_1][\partial'_2]$. Define
$l_i(\phi)=l_\phi(\partial_i)=len(\phi([\partial'_i]))$. If $\phi([\partial'_i])$ is a
hyperbolic isometry in $Isom^+(\H^3)$ then let $\tpartial_i$ denote its axis.

If for $i \ne j$ $\tpartial_i$ and $\tpartial_j$ do not have a common point at infinity
then let $\tpartial_{ij}$ denote the geodesic perpendicular to $\tpartial_i$ and
$\tpartial_j$ oriented from $\tpartial_i$ to $\tpartial_j$. Define
$l_{ij}(\phi)=width(\tpartial_i, \tpartial_j;\tpartial_{ij})$. If orientations on
$\tpartial_i$ and $\tpartial_j$ are not specified, then $l_{ij}(\phi)$ is defined only up
to $\pm i\pi$.

The geodesics $\tpartial_i$ and $\tpartial_{ij}$ for $i,j=1,2,3$ and $i\ne j$ form a
right-angled hexagon. The widths of this hexagon are by definition the numbers $l_i/2$
and $l_{ij}$. It is well known \cite{Fen1989} that up to a rigid motion the hexagon is
determined by the widths $l_1,l_2,l_3$. Therefore the class $[\phi] \in
X(\pi_1(P),PSL_2(\C))$ is determined by $l_1,l_2,l_3$. So it makes sense to define
$l_{ij}(l_1,l_2,l_3)=l_{ij}(\phi)$. When we use this notation, no orientations are
specified, so $l_{ij}$ is well defined only up to $\pm i\pi$.

For later purposes we make the following definitions. If $\phi: \pi_1(P) \to
\Gamma=\pi_1(M)$ is a homomorphism then define $\tv_{ij} \in T_1(\H^3)$ to be the unit
vector based at the point of intersection between $\tpartial_i$ and $\tpartial_{ij}$ and
pointing in the same direction as $\tpartial_{ij}$ (i.e. towards $\tpartial_j$). Let
$v_{ij}(\phi)$ be the projection of $\tv_{ij}$ to $T_1 M$.

By definition $X_{\{P\}}(\pi_1(P),\Gamma)$ is the set of those $[\phi]$ in
$X(\pi_1(P),\Gamma)$ for which $l_i(\phi)$ and $l_{ij}(\phi)$ are defined.

\begin{lem}[Asymptotics for Pants]\label{lem:trig}
For $i=1,2,3$ let $l_i: \R \to \mF$ be continuous functions such that $l_i(t) \to
\infty$ as $t \to \infty$. Let $l_{12}=l_{12}(l_1,l_2,l_3)$. Then
\begin{eqnarray*}
\cosh(l_{12}) &=& 1 + 2e^{-l_1} + 2e^{-l_2} +2e^{(-l_1-l_2+l_3)/2}\\
&&+  O(e^{-2l_1} + e^{-2l_2} + e^{(-l_1-3l_2 + l_3)/2} +
e^{(-3l_1-l_2 + l_3)/2}).
\end{eqnarray*}
\end{lem}

\begin{proof}
This follows from an application of the law of cosines \cite{Fen1989} to the right-angled
hexagon formed from the geodesics $\tpartial_i$ and $\tpartial_{ij}$. Specifically, we
obtain
\begin{eqnarray*}
\cosh(l_{12}) &=& \frac{\cosh(l_3/2) +
\cosh(l_1/2)\cosh(l_2/2)}{\sinh(l_1/2)\sinh(l_2/2)}
\end{eqnarray*}
from which the asymptotics follow.
%Note
%\begin{eqnarray*}
%\frac{1}{\sinh(l_1/2)\sinh(l_2/2)} = 4e^{(-l_1-l_2)/2} +
%O(e^{(-l_1-3l_2)/2} + e^{(-3l_1-l_2)/2})
%\end{eqnarray*}
%and
%\begin{eqnarray*}
%\coth^2(l_1/2)\coth^2(l_2/2) &=& (1+\sinh^{-2}(l_1/2))(1+\sinh^{-2}(l_2/2))\\
%&=& 1 + 4e^{-l_1} + 4e^{-l_2} + O(e^{-2l_1} + O(e^{-2l_2}).
%\end{eqnarray*}
%Hence $\coth(l_1)\coth(l_2) = 1 + 2e^{-l_1 } + 2e^{-l_2} +
%O(e^{-2l_1} + e^{-2l_2})$. So
%\begin{eqnarray*}
%\cosh(l_{12}) &=& 1 + 2e^{-l_1} + 2e^{-l_2} +2e^{(-l_1-l_2+l_3)/2}\\
%&&+  O(e^{-2l_1} + e^{-2l_2} + e^{(-l_1-3l_2 + l_3)/2} +
%e^{(-3l_1-l_2 + l_3)/2}).
%\end{eqnarray*}

\end{proof}

\section{Pulling on Pants Two Legs at a Time}

For hyperbolic elements $g_1, g_2 \in G$ let $X_{\{P\}}(\pi_1(P), G; g_1,g_2)$ be the set
of nondegenerate characters $[\phi]$ such that $\phi([\partial'_i])$ is conjugate to
$g_i$ for $i=1,2$. Similarly, if  $g_1, g_2 \in \Gamma$ let $X_{\{P\}}(\pi_1(P), \Gamma;
g_1,g_2)$ be the set of nondegenerate characters $[\phi] \in X_{\{P\}}(\pi_1(P),\Gamma)$
such that $\phi([\partial'_i])$ is $\Gamma$-conjugate to $g_i$ for $i=1,2$. There is a
projection map $\pi: X_{\{P\}}(\pi_1(P), \Gamma; g_1,g_2) \to X_{\{P\}}(\pi_1(P), G;
g_1,g_2)$.

We put coordinates on $X_{\{P\}}(\pi_1(P), G; g_1,g_2)$ as follows. Fix $i=1,2$. If
$[\phi] \in X_{\{P\}}(\pi_1(P), G; g_1,g_2)$ then by conjugating we may assume that
$\phi([\partial'_i])=g_i$. Let $\tv_{ij}$ be as in the previous section. Recall that
$\tv_{ij} \in T_1^\perp(axis(g_i))$ since $axis(g_i)=\tpartial_i$. But $\tv_{ij}$ is
determined only up to left multiplication by $g_i$. But we will ignore this in order to
simplify the notation. The map
\begin{eqnarray*}
[\phi] \to (\tv_{12}, \tv_{21}, l_3(\phi)) \in T_1^\perp(axis(g_1)) \times
T_1^\perp(axis(g_2)) \times \mF
\end{eqnarray*}
gives coordinates on  $X_{\{P\}}(\pi_1(P), G; g_1,g_2)$.

% Let $vol$ be the measure on $X_{\{P\}}(\pi_1(P), G; g_1,g_2)$ induced by these
%coordinates.

Let $\lambda=\lambda_{g_1,g_2}$ be the measure on $X_{\{P\}}(\pi_1(P), G; g_1,g_2)$ with
derivative
\begin{eqnarray*}
d\lambda &=& \frac{2^{n-3}e^{(n-1)\Re(L-l_1-l_2)/2} }{ (2\pi)^{n-2} vol(T_1 M)}  \,
dvol_{T_1^\perp(axis(g_1))} \, dvol_{T_1^\perp(axis(g_2))} \, dL
\end{eqnarray*}
where $l_i=len(g_i)$.

%For $i=1,2$ let $\sigma_i=\sigma_i(L)$ be a closed oriented geodesic in $M$. Let
%$\lambda=\lambda_L=\lambda_{\sigma_1,\sigma_2}$ be the measure on
%$T_1^\perp(\sigma_1)\times
%%T_1^\perp(\sigma_2) \times \mF$ defined by setting $\lambda(E)$ equal to the number of
%classes $[\phi] \in X_{\{P\}}(\pi_1(P),\Gamma)$ such that
%\begin{itemize}
%\item for $i=1,2$, $\phi[\partial'_i]$ represents $\sigma_i$ and
%\item $\Big(v_{12}(\phi),v_{21}(\phi),l_3(\phi)\Big) \in E$ where $\partial'_i, v_{ij}$
%and $l_3$ are
%defined in section \ref{sec:pants}.
%\end{itemize}

%Let $\lambda'=\lambda'_L$ be the measure on $T_1^\perp(\sigma_1)\times
%T_1^\perp(\sigma_2) \times \mF$ with derivative
%\begin{eqnarray*}
%d\lambda'(E) &=& \frac{2^{n-3}e^{(n-1)\Re(L-l_1-l_2)/2} }{ (2\pi)^{n-2} vol(T_1 M)}  \,
%dvol_{T_1^\perp(\sigma_1)} \, dvol_{T_1^\perp(\sigma_2)} \,
%dL.
%\end{eqnarray*}

\begin{thm}\label{thm:twoleg}
Let $g_1=g_1(t),g_2=g_2(t) \in \Gamma$. If $n=2$ let $\epsilon >0$. If $n=3$, let
$\epsilon \in \A$ be such that $\Re(\epsilon) >0$ and $\Im(\epsilon) >0$. For $i=1,2$ let
$l_i=len(g_i)$ and $L=L(t)\in \mF$. If $L-l_1-l_2 \to \infty$ as $t \to \infty$,
\begin{eqnarray*}
E_L=\sF_1 \times \sF_2 \times B(L,\epsilon) \subset T_1^\perp(axis(g_1))\times
T_1^\perp(axis(g_2)) \times \mF
\end{eqnarray*}
and $\sF_i=\sF'_i \times \sigma_i$, $\sF'_i$ is an arc of $S^{n-2}$ with
$length(\sF'_i)e^{\Re(L-l_1-l_2)/4} \to \infty$ and $\sigma_i$ is a segment of
$axis(g_i)$ with $length(\sigma_i)e^{\Re(L-l_1-l_2)/4} \to \infty$ then $|\pi^{-1}(E_t)|
\sim \lambda_{g_1(t),g_2(t)}(E_t)$ as $t \to \infty$.
\end{thm}

\begin{proof}
By lemma \ref{lem:trig} if $l_1, l_2, l_3, l_3-l_2-l_1$ all tend to infinity then
\begin{eqnarray*}
l_{12}(l_1,l_2,l_3) -\Big( (l_3-l_1-l_2)/2 + \ln(4)\Big) \to 0.
\end{eqnarray*}
By the previous section and theorem \ref{thm:perp} if $n=3$ then
\begin{eqnarray*}
\Big|\pi^{-1}\Big(\sF_1 \times \sF_2 \times B(L,\epsilon)\Big)\Big| &=& \chi\Big(\sF_1
\times \sF_2 \times
l_{12}(l_1, l_2, B(L,\epsilon)\Big)\\
&\sim& \chi'\Big(\sF_1 \times \sF_2 \times B( (L-l_1-l_2)/2 + \ln(4),
\epsilon/2) \Big)\\
%&\sim&  \frac{2^{1-n}vol(\sF_1)vol(\sF_2)}{(2\pi)^{n-2}vol(T_1 M)}\int_{B(
%(L-l_1-l_2)/2 + \ln(4),
%\epsilon/2)  } \,
%e^{(n-1)\Re(L)} \, dL\\
%&\sim& \frac{|\Im(\epsilon)|}{2\pi} \frac{vol(\sF_1)vol(\sF_2)e^{(n-1)\Re(L-l_1-l_2)/2}
%4^{n-1}\sinh((n-1)\Re(\epsilon)/2)}{2^{n-2} (n-1) vol(T_1 M)}\\
&\sim& \frac{2|\Im(\epsilon)|}{(2\pi)^{n-2}} \frac{2^{n-1}\sinh((n-1)\Re(\epsilon)/2)
e^{(n-1)\Re(L-l_1-l_2)/2}vol(\sF_1)vol(\sF_2)}{(n-1) vol(T_1 M)}.
\end{eqnarray*}
A routine calculation shows that this equals $\lambda(E)$. The $n=2$ case is similar.
%But,
%\begin{eqnarray*}
%\lambda(E) &=&  \frac{2^{n-3}}{ vol(T_1 M)}\int_{E_L} \,
%e^{(n-1)\Re(L-l_1-l_2)/2} \, dvol_{T_1^\perp(\sigma_1)} \, dvol_{T_1^\perp(\sigma_2)} \,
%dL\\
%&=&  \frac{2|\Im(\epsilon)|}{(2\pi)^{n-2}}
%\frac{2^{n-1}\sinh((n-1)\Re(\epsilon)/2)e^{(n-1)\Re(L-l_1-l_2)/2}vol(\sF_1)vol(\sF_2)}{
%(n-1) vol(T_1 M)}.
%\end{eqnarray*}

\end{proof}

\section{Pulling on Pants One Leg at a Time}

If $h \in G$ is hyperbolic let $X_{\{P\}}(\pi_1(P), G; h)$ be the set of nondegenerate
characters $[\phi]$ such that $\phi([\partial'_1])$ is conjugate to $h$. Similarly, if
$h \in \Gamma$ let $X_{\{P\}}(\pi_1(P), \Gamma; h)$ be the set of nondegenerate
characters $[\phi] \in X_{\{P\}}(\pi_1(P),\Gamma)$ such that $\phi([\partial'_1])$ is
$\Gamma$-conjugate to $h$. There is a projection map $\pi: X_{\{P\}}(\pi_1(P), \Gamma; h)
\to X_{\{P\}}(\pi_1(P), G; h)$.

We put coordinates on $X_{\{P\}}(\pi_1(P), G; h)$ as follows. Let $\teta=axis(h)$
oriented from the repelling to the attracting fixed point of $h$. Let $\tsigma_2$ be an
oriented geodesic perpendicular to $axis(h)$. The reason for the subscript will be
apparent later. If $[\phi] \in X_{\{P\}}(\pi_1(P), G; h)$ then by conjugating we may
assume that $\phi([\partial'_1])=h$.  Let $\tv_{12}$ be as in section \ref{sec:pants}.
Recall that $\tv_{12} \in T_1^\perp(\teta)$. But $\tv_{12}$ is determined only up to left
multiplication by $h$. Define $y: X_{\{P\}}(\pi_1(P), G; h) \to \mF$ by
$y=width(\tsigma_2, \tv_{12}; \teta)$ where we have abused notation by identifying
$\tv_{12}$ with the oriented geodesic tangent to it.

The map $C:X_{\{P\}}(\pi_1(P), G; h) \to \mF^2 \times \mF/len(h)$ defined by
$C(\phi)=(l_2(\phi),l_3(\phi),y(\phi))$ determines the desired coordinates. To avoid
excessive notation, we identify $X_{\{P\}}(\pi_1(P), G; h)$ with $\mF^2 \times
\mF/len(h)$ through the map $C$. Let $\lambda=\lambda_h$ be the measure on
$X_{\{P\}}(\pi_1(P), G; h)$ with derivative
\begin{eqnarray*}
d\lambda = \frac{1}{(2\pi)^{n-2}}\frac{\exp \Big(\big( \frac{n-1}{2} \big)
\Re(l_2+l_3-len(h))  \Big) }{vol(T_1(M))}  \, dl_2 \, dl_3\, dy.
\end{eqnarray*}

\begin{thm}\label{thm:oneleg}
If $n=2$ let $\epsilon >0$. If $n=3$, let $\epsilon \in \A$ such that $\Re(\epsilon) >0$
and $\pi \ge \Im(\epsilon) >0$. Let $h=h(t) \in \Gamma$, $l_1=len(h)$. Let $l_2,l_3, Y:\R
\to \mF$ be functions of $t$. Let $x=l_{12}(l_1,l_2,l_3)$. Suppose that
\begin{eqnarray*}
E_t= B(l_2,\epsilon) \times B(l_3,\epsilon) \times B(0,Y) \subset \mF^3,
\end{eqnarray*}
$l_2+l_3-l_1 \to \infty$, $\Re(l_3) \ge \Re(l_2)$,
\begin{eqnarray*}
\Big|\Re\big(\sinh(x)\sinh( Y )\big)\Big| e^{\Re(l_2)/2} &\to& \infty \textrm{ and in
case $n=3$ then }\\
\Big|\Im\big(\sinh(x)\sinh( Y )\big)\Big| e^{\Re(l_2)/2} &\to& \infty.
\end{eqnarray*}
Then $|\pi^{-1}(E_t)| \sim \lambda_{h(t)}(E_t)$.
\end{thm}

To prove the theorem, we need to introduce more notation. All quantities and geodesics
below depend implicitly on a parameter $t$. We will derive formulae which are asymptotic
in $t$. Let $[\phi]\in X_{\{P\}}(\pi_1(P),\pi_1(M);h)$.

%We will assume throughout that $\phi([\partial'_1])$ represesents $\eta$ and $\teta :=
%axis \, \phi([\partial'_1])$ is a lift of $\eta$ to $\H^n$.

%Let $\sigma$ be a geodesic segment perpendicular to $\eta$. Let $\teta$ be a lift of
%$\eta$.

\begin{figure}[htb]
\begin{center}
\ \psfig{file=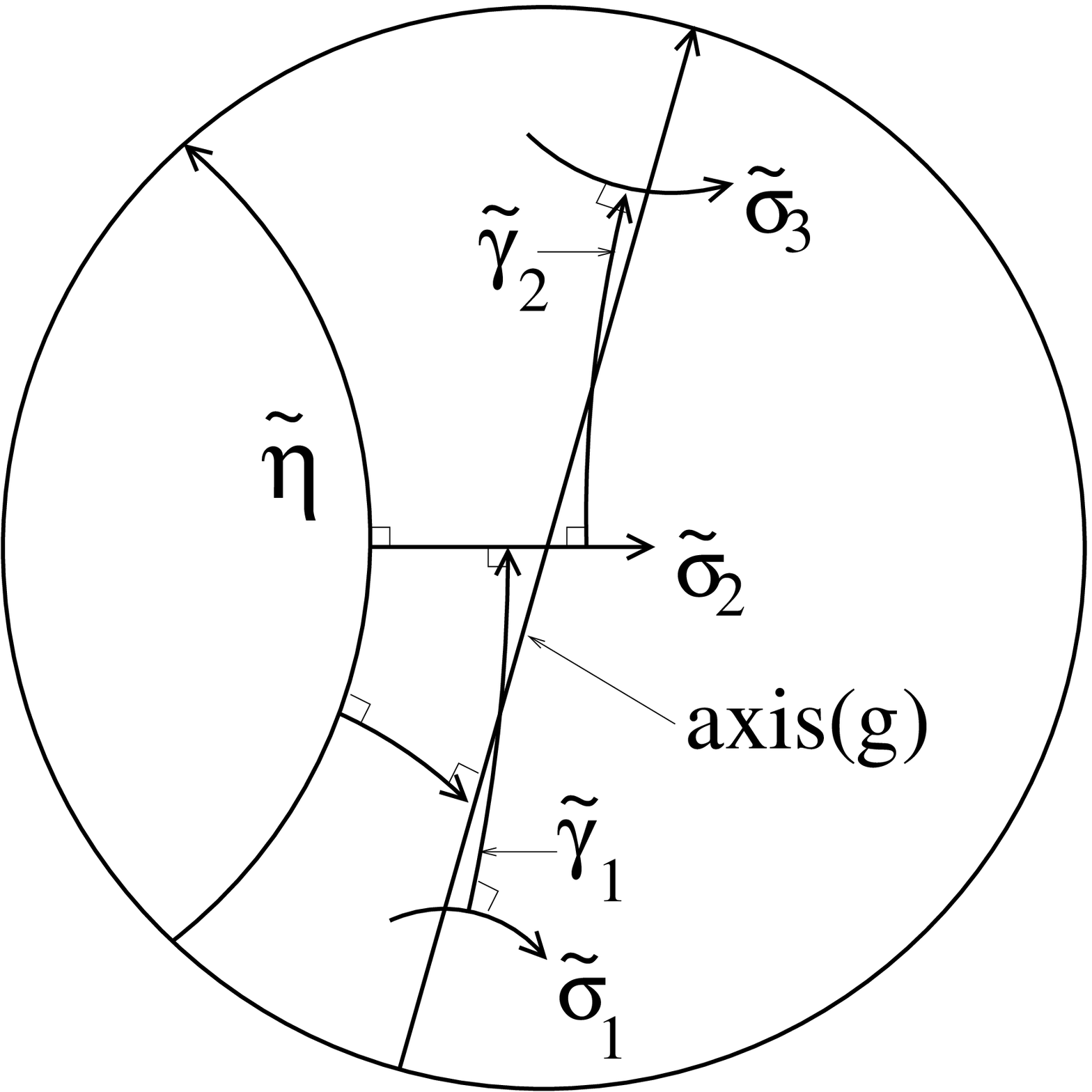,height=2in,width=2in}
\caption{}
\label{fig:stairs}
\end{center}
\end{figure}

Let $g = \phi \, [\partial'_2]$. Let $\tsigma_1 = g^{-1}\tsigma_2$ and
$\tsigma_3=g\tsigma_2$. See figure \ref{fig:stairs}. For $i=1,2$ let $\tgamma_i$ be the
segment perpendicular to $\tsigma_i$ and $\tsigma_{i+1}$. Let $\tpartial_{12}$ be
perpendicular to $\teta$ and $axis(g)$. Let $\tau$ be perpendicular to $\tsigma_2$ and
$axis(g)$. Define
\begin{displaymath}
\begin{array}{lll}
x=l_{12}(\phi)=width(\eta,axis(g);\tpartial_{12}) &
y=width(\tsigma_2,\tpartial_{12};\teta) & z=width(\tsigma_2,axis(g);\tau)\\
u_1=width(\teta,\tgamma_1;\tsigma_2) & u_2=width(\teta, \tgamma_2;\tsigma_2)
&l'_2=len(\gamma)= width(\tsigma_1,\tsigma_2;\tgamma_1)\\
m=\frac{u_1+u_2}{2} & w=u_2-u_1 & l_i=l_i(\phi).
\end{array}
\end{displaymath}

\begin{lem}\label{lem:cov}
If $l_1,l_2,l_3 \to \infty$, $l_2+l_3-l_1 \to +\infty$, $\Re(l_3) \ge \Re(l_2)$ and
$\tanh(x)\sinh(y) \to 0$ as $t \to \infty$ then
\begin{eqnarray*}
l'_2 -l_2 &\to& 0\\
m &\sim& x\sim \arccosh(1 + 2e^{(l_3-l_1-l_2)/2}) \\
w &\sim&  -2\sinh(m)\sinh(y).
\end{eqnarray*}
\end{lem}
After this lemma is proven, the theorem will follow from theorem \ref{thm:perp} and a
change of variables. We need some intermediate lemmas.

%Note $w=u_2-u_1=width(\tgamma_1,\tgamma_2;\tsigma_2)$.

\begin{lem}\label{lem:cov1}
If $w\to 0$ and $l_2 \to \infty$ (as $t \to \infty$) then $len(g) - len(\gamma) \to  0$
(i.e. $l_2-l'_2 \to 0$) and $z -i\pi/2\sim -iw/2$.
\end{lem}

\begin{proof} Recall that a right-angled hexagon is an ordered $6$-tuplet of oriented
geodesics $(\tH_1,..,\tH_6)$ such that for all $i$ mod $6$, $\tH_i$ is perpendicular to
$\tH_{i+1}$ \cite{Fen1989}. Let $H=(\tH_1,...,\tH_6)$ be the right-angled hexagon
satisfying
\begin{displaymath}
\begin{array}{llllll}
\tH_1= \tgamma_1, & \tH_2= \tsigma_2, & \tH_3 = \tau, & \tH_4 = axis(g), & \tH_5 =
-(g^{-1} \tau), & \tH_6 = \tsigma_1.
\end{array}
\end{displaymath}
\begin{figure}[htb]
\begin{center}
\ \psfig{file=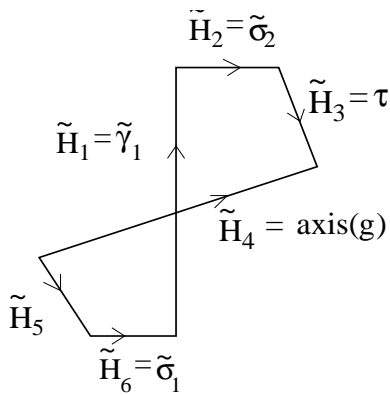,height=2in,width=2in}
\caption{A nongeometric diagram of hexagon $H$.}
\label{fig:hexH}
\end{center}
\end{figure}
See figure \ref{fig:hexH}. Let $H_i:=width(\tH_i,\tH_{i+2};\tH_{i+1})$ (subscripts mod
6). We claim that
\begin{displaymath}
\begin{array}{lll}
H_1 = len(\gamma) & H_2 = w/2 +i\pi/2 & H_3 = z\\
H_4 = -len(g) + i\pi & H_5 = z & H_6 = w/2 + i\pi/2.
\end{array}
\end{displaymath}
The values of $H_1$ and $H_3$ follow by definition. $g$ maps
$\tH_4$ to itself, $\tH_5$ to $-\tH_3$ and $\tH_6$ to $\tH_2$. Since
\begin{eqnarray*}
width(gH_i,gH_j;gH_k)=width(H_i,H_j;H_k)
\end{eqnarray*}
for any $i,j,k$ it follows that
\begin{eqnarray*}
H_5=width(\tH_4,\tH_6;\tH_5) = width(\tH_4,
\tH_2; -\tH_3) = width(\tH_2,\tH_4;\tH_3) = H_3 = z.
\end{eqnarray*}
Since $g$ maps $\tH_5$ to $-\tH_3$ and $axis(g)=\tH_4$,
\begin{eqnarray*}
len(g)=width(\tH_5,-\tH_3;\tH_4)= width(\tH_5,\tH_3;\tH_4)+i\pi = -H_4 +i\pi.
\end{eqnarray*}
Hence $H_4=-len(g)+i\pi$. Observe that
\begin{eqnarray*}
w&=&u_2-u_1=width(\tgamma_1,\tgamma_2;\tsigma_2) =
width(\tgamma_1,\tau;\tsigma_2)+width(\tau,\tgamma_2;\tsigma_2)\\
&=&width(\tH_1,\tH_3;\tH_2) + width(g^{-1} \tau,g^{-1}\tgamma_2;g^{-1}\tsigma_2)\\
&=&H_2 + width(-\tH_5, \tH_1; \tH_6) = H_2 + H_6 + i\pi.
\end{eqnarray*}
The law of sines implies
\begin{eqnarray*}
\frac{\sinh(H_3)}{\sinh(H_6)} = \frac{\sinh(H_5)}{\sinh(H_2)}.
\end{eqnarray*}
Since $H_3=H_5$ this implies $\sinh(H_6)=\sinh(H_2)$. Now we can use the law of
cosines to obtain $\cosh(H_6)=\cosh(H_2)$. Together these imply
$H_2=H_6$. Thus, $H_2 =H_6 = w/2 \pm i\pi/2$. From the figure, it can be
checked that $H_2=H_6=w/2 + i\pi/2$. This proves the claim.

% this accords with Fenchel's
%weird convention. we should explicate this in the background
%section. It's not so weird after all: a unit normal to a plane induces
%an orientation on the plane in the obvious way. If $H_{i+1}$ is normal
%to the plane spanned by $H_i$ and $H_{i+2}$ then we use this
%orientation to measure the angle from $H_i$ to $H_{i+2}$. I had this
%backwards before.

By the law of cosines,
\begin{eqnarray*}
\cosh(H_4) &=& \cosh(H_2)\cosh(H_6) + \sinh(H_2)\sinh(H_6)\cosh(H_1).
\end{eqnarray*}
Substitution yields
\begin{eqnarray*}
\cosh(-len(g)+i\pi) &=& \cosh^2(w/2+i\pi/2) + \sinh^2(w/2+i\pi/2)\cosh(len(\gamma)).
\end{eqnarray*}
This simplifies to
%\begin{eqnarray*}
%-\cosh(len(g)) &=& -\sinh^2(w/2) - \cosh^2(w/2)\cosh(len(\gamma)).
%\end{eqnarray*}
%And,
\begin{eqnarray*}
\cosh(len(g)) &=& \cosh(len(\gamma)) + \sinh^2(w/2)(1 +
\cosh(len(\gamma)))\\
&=& \cosh(len(\gamma)) + O(w^2e^{len(\gamma)}).
\end{eqnarray*}
This implies $len(g)=len(\gamma)+O(w^2)$. By the law of cosines
\begin{eqnarray*}
\cosh(H_3) &=& \cosh(H_1)\cosh(H_5) + \sinh(H_1)\sinh(H_5)\cosh(H_6).
\end{eqnarray*}
Substitution yields
\begin{eqnarray*}
\cosh(z) &=& \cosh(len(\gamma))\cosh(z) +
               \sinh(len(\gamma))\sinh(z)\cosh(w/2 + i\pi/2).
\end{eqnarray*}
By solving for $\coth(z)$ we obtain
\begin{eqnarray}\label{eqn:cothz}
\coth(z) &=& \frac{\sinh(len(\gamma))\cosh(w/2 +
i\pi/2)}{1-\cosh(len(\gamma))}.
\end{eqnarray}
Since $l_2=len(g) \to \infty$, $len(\gamma)\to \infty$. So $z-i\pi/2 \sim -iw/2$.
%
%Now, $\sinh(len(\gamma))/(1-\cosh(len(\gamma))) = -1 + O(w^2)$
%and $\sinh(z)=i + O(w^2)$. Hence
%\begin{eqnarray*}
%\cosh(z) &=& -i \cosh(w/2 + i\pi/2) + O(w^2)\\
%        &=& \sinh(w/2) + O(w^2)\\
%&=& w/2 + O(w^2)
%\end{eqnarray*}
%But $z=i\pi/2 + O(w^2)$. Hence
%\begin{eqnarray*}
%z -i\pi/2 &=& \sinh(z - i\pi/2) + O(w^2)\\
%&=& -i \cosh(z) + O(w^2)\\
%&=& -iw/2 +  O(w^2).
%\end{eqnarray*}
%So $z= -iw/2 + i\pi/2 + O(w^2)$.

%By the law of sines:
%\begin{eqnarray*}
%\frac{\sinh(H_1)}{\sinh(H_4)} = \frac{\sinh(H_3)}{\sinh(H_6)}
%\end{eqnarray*}
%Hence
%\begin{eqnarray*}
%\frac{\sinh(len(\gamma))}{\sinh(-len(g)+i\pi)} =
%\frac{\sinh(z)}{\sinh(w/2+i\pi/2)}.
%\end{eqnarray*}

%So:
%\begin{eqnarray}\label{eqn:sines}
%\frac{\sinh(len(\gamma))}{\sinh(len(g))} =
%\frac{\sinh(z)}{i\cosh(w/2)}.
%\end{eqnarray}

%Since $w=O(1)$ it follows that
%$\sinh(w/2+i\pi/2)=i\cosh(w/2) = i + O(w^2)$. It follows from the
%above and $len(g)=len(\gamma)+O(w^2)$ that $\sinh(z)=i +
%O(w^2)$. Therefore $z=i\pi/2 + O(w^2)$.

\end{proof}

\begin{lem}\label{lem:cov2}
If $\tanh(x)\sinh(y) \to 0$ and $l_2 \to \infty$ then $w \to 0$, $x \sim m$ and $w \sim
-2 \sinh(m)\sinh(y)$.
\end{lem}

\begin{proof}
Let $P=(\tP_1,..,\tP_5)$ be the right-angled pentagon defined by
\begin{displaymath}
\begin{array}{lllll}
\tP_1= \tsigma_2, & \tP_2= \tau, & \tP_3 = axis(g), & \tP_4 = \tpartial_{12}, & \tP_5 =
\teta.
\end{array}
\end{displaymath}
\begin{figure}[htb]
\begin{center}
\ \psfig{file=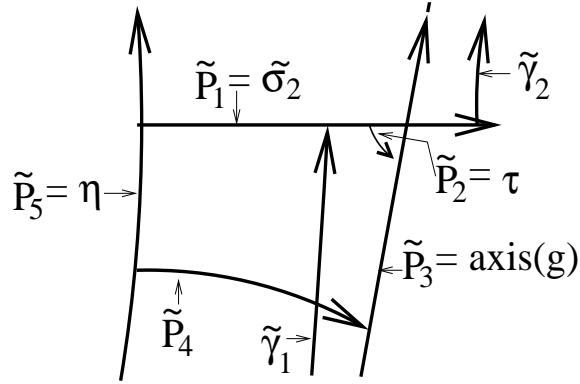,height=2in,width=3in}
\caption{A diagram of pentagon $P$.}
\label{fig:pentagon}
\end{center}
\end{figure}
See figure \ref{fig:pentagon}. Let $P_i:=(\tP_{i-1},\tP_{i+1};\tP_i)$ (subscripts mod 5).
We claim
\begin{displaymath}
\begin{array}{llll}
P_1 = m + i\pi/2, & P_2 = z, & P_4 = -x, & P_5 = -y.
\end{array}
\end{displaymath}
The last three identities hold by definition. The first follows from
\begin{eqnarray*}
P_1&=&width(\teta,\tau;\tsigma_2) = width(\teta, \tgamma_1;\tsigma_2) +
width(\tgamma_1,\tau;\tsigma_2)\\
&=& u_1 + H_2 = u_1 + w/2 + i\pi/2 = u_1 + \frac{u_2-u_1}{2} + i\pi/2 = m + i\pi/2.
\end{eqnarray*}
By the trigonometric formulae for right angled pentagons \cite{Fen1989}
\begin{eqnarray*}
\cosh(P_4) &=& -\sinh(P_1)\sinh(P_2).
\end{eqnarray*}
Substituting yields $\cosh(-x) = -\sinh(m+i\pi/2)\sinh(z)$. Equivalently,
\begin{eqnarray}\label{eqn:sinhz}
\sinh(z)&=& \frac{i \cosh(x)}{\cosh(m)}.
\end{eqnarray}
By other trigonometric formulae for right-angled pentagons \cite{Fen1989},
\begin{eqnarray}\label{eqn:coshz}
\cosh(z) &=& \cosh(P_2) = -\sinh(P_5)\sinh(P_4)\\
&=& -\sinh(-y)\sinh(-x) = - \sinh(x)\sinh(y).
\end{eqnarray}
Dividing equation \ref{eqn:coshz} by equation \ref{eqn:sinhz} yields $\coth(z) = i
\tanh(x)\sinh(y)\cosh(m)$. With equation \ref{eqn:cothz} this implies
\begin{eqnarray*}
i\tanh(x)\sinh(y) &=& -\frac{\sinh(len(\gamma))\sech(m)\cosh(w/2 +
i\pi/2)}{\cosh(len(\gamma))-1)}.
\end{eqnarray*}
By hypothesis $\tanh(x)\sinh(y)\to 0$. Since
$\frac{\sinh(len(\gamma))}{\cosh(len(\gamma))-1}$ and $\sech(m)$ are bounded from below,
this implies that $\sinh(w) \to 0$. So $w\to 0$. This allows us to use results from the
previous lemma.

Since $\sinh(z) \to i$ and $\sinh(m+i\pi/2) = i\cosh(m)$ equation \ref{eqn:sinhz} implies
$\cosh(x) -1 \sim \cosh(m)-1$. This implies $x\sim \pm m$. From the figure it can be seen
that $x \sim m$. By equation \ref{eqn:coshz}
\begin{eqnarray*}
\sinh(y) &=& \frac{-\cosh(z)}{\sinh(x)} \sim  \frac{- \cosh(-iw/2 + i\pi/2)}{\sinh(m)}=
\frac{ i \sinh(iw/2) }{\sinh(m)} \sim \frac{-w}{2\sinh(m)}.
\end{eqnarray*}
So $w \sim -2\sinh(m)\sinh(y)$.

\end{proof}

\begin{proof}(of lemma \ref{lem:cov})
By lemma \ref{lem:cov2} $w \to 0$ so we can use lemma \ref{lem:cov1}. So $l'_2-l_2 \to
0$. Since $l_2+l_3 - l_1 \to \infty$, $\Re(l_3) \ge \Re(l_2)$ and $l_1,l_2,l_3 \to
\infty$ it follows from lemma \ref{lem:trig} that
\begin{eqnarray*}
\cosh(l_{12}(\phi))-1 \sim 2e^{(l_3-l_1-l_2)/2}.
\end{eqnarray*}
By definition, $x=l_{12}(\phi)$. By lemma \ref{lem:cov2} $x\sim m$ so $\cosh(m)-1 \sim
2e^{(l_3-l_1-l_2)/2}$ too. By lemma \ref{lem:cov2} $w \sim -2\sinh(m)\sinh(y)$.
\end{proof}

\begin{proof}(of theorem \ref{thm:oneleg})

Let $\eta$ be the projection of $\teta$ to $M=\H^n/\Gamma$. Let $\sigma$ be the
projection of $\tsigma_2$ to $M$.

Let $\chi=\chi_{\sigma,\sigma}$ be as defined in theorem \ref{thm:perp}. Recall that
$\chi$ is a measure on
$T_1^\perp(\sigma) \times T_1^\perp(\sigma) \times \mF$. Indeed $\chi(E)$ is the number of
oriented perpendiculars $\gamma$ from $\sigma$ to $\sigma$ such that
$(v_1,v_2,len(\gamma)) \in E$ where $v_i$ is the $i$-th endpoint vector of $\gamma$. It
is convenient for the proof to put coordinates on $T_1^\perp(\sigma)$. So define $u:
T_1^\perp(\sigma) \to \mF$ by $u(v)=width(\eta, v; \sigma)$. We have abused notation here
by identifying $v$ with the oriented geodesic tangent to it. Define $K:
T_1^\perp(\sigma)^2 \times \mF \to \mF^3$ by $K(v_1,v_2,l_2) = \big(l_2,
u(v_1),u(v_2)\big)$.

Let $\rho$ be the density of $K_*\chi'$. From the definition of $\chi'$ (see the
paragraphs before theorem \ref{thm:perp}) it follows that
\begin{eqnarray*}
\rho(L, u_1,u_2) &=&   \frac{e^{(n-1)\Re(L)}}{(2\pi)^{n-2} 2^{n-1}vol(T_1 M)}.
\end{eqnarray*}
%By theorem \ref{thm:perp}, $(K_*\chi)(E_L) $ is asymptotic to $(K_*\chi')(E_L)$ when
%$\{E_L\}$
%satisfies hypotheses corresponding to those of theorem \ref{thm:perp}.

The parameters $l_1, l_2, l_3, y$ determine $\teta, \tsigma_2, axis(g), \tgamma_i$, etc.
up to a rigid motion. Hence, keeping $l_1$ fixed we may define $G(l_2,l_3,y)=(l_2', m,
w)$. Define
\begin{eqnarray*}
H(l_2',m,w)=(l'_2,u_1,u_2) = \Big(l'_2, m-\frac{w}{2}, m+ \frac{w}{2}\Big).
\end{eqnarray*}
By definition $\lambda(E_t) = (K_* \chi)(H \circ G (E_t))$. It follows from theorem
\ref{thm:perp} that  $ (K_* \chi)(H \circ G (E_t)) \sim (K_* \chi')(H \circ G (E_t))$
once we show that its hypotheses are satisfied. It is immediate that $K^{-1}(H \circ
G(E_t))$ tends to infinity with $t$. We need to show that
\begin{eqnarray*}
\lim_{t \to \infty} \frac{\chi'(N_1(1, \partial K^{-1} \circ H \circ G(E_t)))}{\chi'(
K^{-1} \circ H \circ G(E_t)) )} =0.
\end{eqnarray*}
By definition
\begin{eqnarray*}
N_1(1, \partial K^{-1} \circ H \circ G(E_t))&:=& K^{-1}\Big(\{ (u_1,u_2,l_2) \in \mF^3 ;
\, \exists (u'_1,u'_2,l_2) \in \partial (H \circ G)(E_t) \\
&& \textrm{ such that } d(u'_1,u_1) \le  e^{-\Re(l_2)/2}, d(u'_2,u_2) \le
e^{-\Re(l_2/2)}\}\Big).
\end{eqnarray*}
But the Jacobian of the map $(u_1,u_2) \to (m,w) =(\frac{u_1+u_2}{2},u_2-u_1)$ is
constant. So it suffices to show that
\begin{eqnarray}\label{eqn:N1}
\lim_{t \to \infty} \frac{\chi'\Big(K^{-1} \circ H \big(N_1(1, \partial G(E_t) )\big)
\Big)}{\chi'( K^{-1} \circ H \circ G(E_t) )} =0.
\end{eqnarray}
where
\begin{eqnarray*}
N_1(1, \partial G(E_t))&:=& \{ (l_2,m,w) \in \mF^3 ; \, \exists (l_2,m',w') \in \partial
G(E_t) \\
&& \textrm{ such that } d(m,m') \le  e^{-\Re(l_2)/2}, d(w,w') \le e^{-\Re(l_2/2)}\}.
\end{eqnarray*}
By lemma \ref{lem:cov}, $m \sim \arccosh(1+2\exp( (l_3-l_1-l_2)/2))$. Since $l_3+l_2-l_1
\to \infty$ it follows that the range of $\Re(m) e^{\Re(l_2)/2} \to \infty$ and, if
$n=3$, the range of $\Im(m)e^{\Re(l_2)/2} \to \infty$. Since $w \sim -2\sinh(m)\sinh(y)$,
and $\Re(\sinh(m)\sinh(Y))e^{\Re(l_2)/2} \to \infty$ (and if $n=3$,
$\Im(\sinh(m)\sinh(Y))e^{\Re(l_2)/2} \to \infty$ too) by hypothesis, equation
\ref{eqn:N1} is satisfied.

The last hypothesis of theorem \ref{thm:perp} is easy to verify. So $(K_*\chi)(H \circ
G(E_L)) \sim (K_*\chi')(H \circ G(E_L))$. Define
\begin{eqnarray*}
\tG(l_2,l_3,y) &=& \Big(l_2, \arccosh(1+2e^{(l_3-l_1-l_2)/2}),
-2\sinh\big(\arccosh(1+2e^{(l_3-l_1-l_2)/2})\big) \sinh(y) \Big).
\end{eqnarray*}
Lemma \ref{lem:cov} implies $\tG \sim G$ provided that $\tanh(x)\sinh(y) \to 0$ for $y
\in B(0,Y)$. This may not be the case, but we can always choose $Y$ smaller if necessary
so that it holds. But we cannot choose $Y$ too small since it is necessary to maintain
the hypothesis that
\begin{eqnarray*}
\Re(\sinh(x)\sinh( Y ))e^{\Re(l_2)/2} \to \infty \textrm{ and if $n=3$ then }
\Im(\sinh(x)\sinh( Y ))e^{\Re(l_2)/2} \to \infty \textrm{, too.}
\end{eqnarray*}
It is easy to see that this can be achieved.

Now that we have $G \sim \tG$ it can easily be verified that $(K_*\chi')(H \circ G(E_L))
\sim (K_*\chi')(H \circ \tG(E_L))$. So we have shown that $|\pi^{-1}(E_t)| \sim
(K_*\chi')(H \circ \tG(E_L))$. To finish we need only compute the later:
\begin{eqnarray*}
(K_*\chi')(H \circ \tG(E_t)) &=& \int_{H\circ G(E_t)} \rho(l_2, u_1, u_2) dl_2 \, du_1 \,
du_2\\
   &\sim& \int_{E_t} \rho(H \circ \tG(l_2,l_3,y)) |Jac(H \circ \tG)| dl_2 \, dl_3 \, dy.
\end{eqnarray*}
An easy computation shows that the Jacobian of $H$ is $1$ and
\begin{eqnarray*}
|Jac(\tG)| &=& 2^{n-1} e^{(n-1)\Re(l_3-l_1-l_2)/2}.
\end{eqnarray*}
Another short computation shows
\begin{eqnarray*}
\rho(H \circ \tG(l_2,l_3,y)) &=& \frac{ e^{(n-1)\Re(l_2)}}{(2\pi)^{n-2}2^{n-1} vol(T_1
M)}.
\end{eqnarray*}
The theorem follows from the above three equations.

\end{proof}

\section{Telescoping Paths}\label{sec:telescope}

\begin{thm}\label{thm:telescope}
If $F$ is a free group of finite rank and for each $t>0$, $[\phi_t] \in X_{\{P\}}(F, Isom^+(\H^2))$ is a
sequence of discrete and faithful characters such that the length of the shortest closed
curve on $\H^2/\phi_t(F)$ tends to infinity then there exists a pants decomposition
$\sP_t$ of the convex core of $\H^2/\phi_t(F)$ such that $\{([\phi_t],\sP_t)\}_{t>0}$ is
telescoping.
\end{thm}

It is possible that this result holds in all dimensions.

\begin{proof}
Let $S_t$ be the convex core of $\H^2/\phi_t(F)$. By definition, the
convex core is the smallest closed convex subset of $\H^2/\phi_t(F)$ that is
homotopy equivalent to it. It is a compact surface with geodesic boundary.

The proof is by induction on the number of pants in a pants
decomposition of $S_t$. If $S_t$ is a single pair of pants, then the
result is vacuous.

Let $\alpha_t$ be the shortest homotopically nontrivial arc on $S_t$ with both
endpoints in the boundary. Let $P_t \subset S_t$ be the unique
pair of pants containing $\alpha_t$. If $b_{1}, b_{2}$ are the boundary components
containing the endpoints of $\alpha$ then $P_t$ is homotopy equivalent to a regular
neighborhood of $\alpha_t \cup b_{1} \cup b_{2}$. It is possible that $b_{1}=b_{2}$.

Let $S'_t$ be the closure of
$S_t-P_t$. Let $F'$ be the free group with $rank(F')=rank(F)-1$. There exists a
representation $\phi'_t: F' \to Isom^+(\H^2)$ such that $S'_t$ is isometric to the convex
hull of
$\H^2/\phi'_t(F')$. By induction, we may assume that there exists a pants decomposition
$\sP'_t$ of $S'_t$ for which $(\phi'_t, \sP'_t)$ is telescoping. So it suffices to show
that $len(\partial S_t) - len(\partial S'_t) \to \infty$ as $t \to \infty$.

We claim that the length of $\alpha_t$ tends to zero as $t\to
\infty$. Since the length of the shortest closed geodesic on $S_t$ tends to infinity it
follows that
$len(\partial S_t) \to \infty$. Since $\alpha_t$ is the shortest
nontrivial arc with endpoints in the boundary, it follows that the
$len(\alpha_t)/2$ neighborhood of the boundary does not have
self-intersections. Thus its area is at least $len(\partial S_t)len(\alpha_t)/2$. But the
total area of the surface is bounded by $2 \pi (rank(F)-1)$ by Gauss-Bonet. So
$len(\partial S_t)len(\alpha_t)/2 \le 2 \pi (rank(F)-1)$. Since $len(\partial S_t)$ tends
to infinity, it must be that $len(\alpha_t) \to 0$.

There are two cases depending on whether $P_t \cap S'_t$ has one or
two components.

{\bf Case 1}: Suppose that $P_t \cap S'_t$ is a single component which we
call $\partial_3$. Then $\alpha_t$ must be the
shortest arc between the other two boundary components of $P_t$. In the
notation of lemma \ref{lem:trig}, the other two components are
$\partial_1, \partial_2$ with lengths $l_1(\phi_t), l_2(\phi_t)$ and
$\alpha_t=\partial_{12}$ is the shortest arc from $\partial_1$ to $\partial_2$. Thus
$len(\alpha_t)=l_{12}$. Thus we have
\begin{eqnarray*}
\cosh(len(\alpha_t)) &=& 1 + 2e^{-l_1} + 2e^{-l_2} +2e^{(-l_1-l_2+l_3)/2}\\
&&+  O(e^{-2l_1} + e^{-2l_2} + e^{(-l_1-3l_2 + l_3)/2} +
e^{(-3l_1-l_2 + l_3)/2})
\end{eqnarray*}
where $l_3$ is the length of $\partial_3$. Since $len(\alpha_t)$ tends
to zero, it must be that $l_1+l_2-l_3$ tends to infinity. But
$l_1+l_2-l_3 = len(\partial S_t) - len(\partial S'_t)$. This finishes case 1.

{\bf Case 2}: Suppose that $P_t \cap S'_t$ has two components
$\partial_1$ and $\partial_2$. Let $\partial_3$ be the other boundary component. Both
endpoints of $\alpha_t$ lie in $\partial_3$. As in section
\ref{sec:pants}, let $\partial_{ij}$ be the shortest segment from $\partial_i$ to
$\partial_j$. $\alpha_t$ is the union of the two shortest paths from $\partial_3$ to
$\partial_{12}$. Let $l_{ij}$ be the length of $\partial_{ij}$.

Cut the pair of pants $P$ along the segments $\partial_{ij}$ for all $i\ne j$. This
decomposes $P$ into two isometric right-angled hexagons in the usual
way. Next, cut along $\alpha_t$. This decomposes each hexagon into 2 right-angled
pentagons one of which has 3 different sides of length $l_1/2$, $l_{13}$ and
$len(\alpha_t)/2$. Let $l_\alpha:=len(\alpha_t)$. See figure \ref{fig:telescope}.
\begin{figure}[htb]
\begin{center}
\ \psfig{file=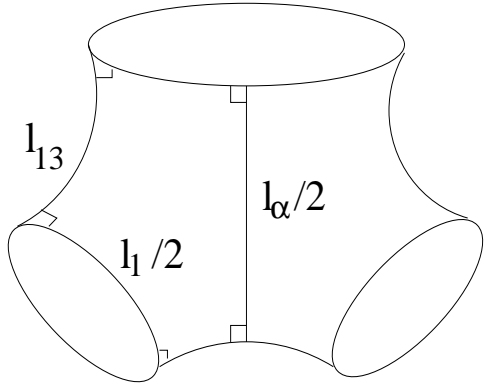,height=2in,width=2.5in}
\caption{The pairs of pants $P$ with pentagon shown and 3 sides labeled by their lengths.}
\label{fig:telescope}
\end{center}
\end{figure}

From the trigonometric formulas for right-angled pentagons \cite{Fen1989}, we have
\begin{eqnarray}\label{eqn:a}
\cosh(l_\alpha/2) &=& \sinh(l_1/2)\sinh(l_{13}).
\end{eqnarray}
%Since $l_\alpha \to 0$, $\cosh(l_\alpha/2) \to 1$. Since $l_1$ is the length of a closed
%geodesic on $S_t$, the hypotheses imply $l_1 \to \infty$. So it must be that
%$\sinh(l_{13})$, and therefore $l_{13}$, tend to
%zero. So $\cosh(l_{13})$ tends to $1$.
By the law of cosines
\begin{eqnarray*}
\cosh(l_{13}) &=& \frac{\cosh(l_2/2)}{\sinh(l_1/2)\sinh(l_3/2)} +
\coth(l_1/2)\coth(l_3/2).
\end{eqnarray*}
%Since $l_1$ and $l_3$ both tend to infinity, $\coth(l_1/2)$ and $\coth(l_3/2)$ tend to
%$1$.
Square both sides of equation \ref{eqn:a} to obtain
\begin{eqnarray}\label{eqn:a2}
\cosh^2(l_\alpha/2) &=& \sinh^2(l_1/2)(\cosh^2(l_{13})-1)\\
                &=& \frac{ \sinh^2(l_1/2)\cosh^2(l_2/2)}{\sinh^2(l_1/2)\sinh^2(l_3/2)}\\
&&+ 2\frac{\sinh^2(l_1/2)\cosh(l_2/2)\coth(l_1/2)\coth(l_3/2)}{\sinh(l_1/2)\sinh(l_3/2)}\\
&& + (\coth^2(l_1/2)\coth^2(l_3/2)-1)\sinh^2(l_1/2).
\end{eqnarray}
The last term above equals
\begin{eqnarray*}
\cosh^2(l_1/2)\coth^2(l_3/2)-\sinh^2(l_1/2) &=& \sinh^2(l_1/2)(\coth^2(l_3/2)-1) +
\coth^2(l_3/2) \ge 1.
\end{eqnarray*}
Since $\cosh^2(l_\alpha/2) \to 1$ the other terms in equation
\ref{eqn:a2} tend to zero. In particular,
\begin{eqnarray*}
2\frac{\sinh^2(l_1/2)\cosh(l_2/2)\coth(l_1/2)\coth(l_3/2)}{\sinh(l_1/2)\sinh(l_3/2)} \to
0.
\end{eqnarray*}
This implies $l_3 - l_1 - l_2$ tends to infinity. But $l_3-l_1-l_2 =
len(\partial S_t) - len(\partial S'_t)$. This completes case 2 and finishes the theorem.

\end{proof}

\section{Proofs of the Main Theorems}

\begin{proof}(of theorem \ref{thm:3d})
As in the proof of lemma \ref{lem:reduction} we may assume without loss of generality
that $E_t \subset X_\sP(\pi_1(S),G)$ is a rectangle with respect to the Fenchel-Nielsen
coordinates induced by $\sP$.

The proof is by induction on the rank $r$ of the free group $\pi_1(S)$. Recall
$\sP=\{P_{i}\}_{i=1}^{r-1}$. If $r=2$ let $S'=S'_t$ be a boundary component of $P_1$ such
that $\Re(l_{\phi_t}(S')) \le \Re( l_{\phi_t} (\gamma))$ for any other boundary component
$\gamma$ of $P_1$. Otherwise let $S'=\cup_{i=1}^{r-2} P_i$.

Define the restriction map
\begin{eqnarray*}
R: X_\sP(\pi_1(S), G) \to X_{\sP'}(\pi_1(S'),G) &,& R([\phi]) = [\phi|_{\pi_1(S')}]
\end{eqnarray*}
where $\sP'=\{P_i\}_{i=1}^{r-2}$ if $r>2$ and $\sP'=\emptyset$ otherwise. If $r=2$,
$X_{\sP'}(\pi_1(S'),G)$ is just $X(\pi_1(S'),G)$. If $r=2$ let $\nu'$ be the measure on
$X_{\sP'}(\pi_1(S'),G)$ with derivative
\begin{eqnarray*}
d\nu' = \frac{e^{(n-1)\Re( len(S') )}}{(2\pi)^{n-2}\Re( len(S') )} \, dlen(S').
\end{eqnarray*}
If $r>2$ then let $\nu'$ be the measure on $X_{\sP'}(\pi_1(S'),G)$ with derivative
\begin{eqnarray*}
d\nu' = vol(T_1 M)^{2-r} 2^{(n-3)(-genus(S'))} (2\pi)^{(1-r)(n-2)}
e^{\Re[len_{\psi}(\partial S)](n-1)/2} \, dvol_{\sP'}(\psi).
\end{eqnarray*}
We claim $|\pi^{-1}(R(E_t))| \sim \nu'(R(E_t))$. This follows from the prime geodesic
theorem if $r=2$ and from the induction hypothesis if $r>2$.

To simplify notation, let $P$ denote $P_{r-1}$. Suppose $[\phi'] \in
X_{\sP'}(\pi_1(S'),G)$ and $[\phi''] \in X_{\{P\}}(\pi_1(P), G)$ are such that for every
$g \in \pi_1(S)$ which represents a curve in $P \cap S'$, $\phi'(g)$ is conjugate to
$\phi''(g)$. Then we may form the connected sum
\begin{eqnarray*}
\phi=\phi' \#_{P \cap S'} \, \phi'' \in X_\sP(\pi_1(S),G)
\end{eqnarray*}
in the obvious way. For $\phi' \in R(E_t)$, let $E_{t,\phi'}$  be the set of all $\phi''
\in X_{\{P\}}(\pi_1(P), G)$ for which $\phi' \#_{P \cap S} \, \phi'' \in E_t$. So
$E_{t,\phi'} = R^{-1}(\phi') \cap E_t$.

Let $c$ be the number of components of $P \cap S$. If $c=1$ then $P \cap S$ has a single
component $\partial_1$ and $E_{t,\phi'} \subset X_{\{P\}}(\pi_1(P),G;\phi'(\partial_1))$.
If $c=2$ then $P \cap S$ consists of two components $\partial_1,\partial_2$ and
$E_{t,\phi'} \subset X_{\{P\}}(\pi_1(P),G;\phi'(\partial_1),\phi'(\partial_2) )$. We have
abused notation here by identifying $\partial_1$ and $\partial_2$ with elements of
$\pi_1(P)$ that represent them.

If $c=1$, let $\lambda=\lambda_{\phi'(\partial_1)}$ be as in theorem \ref{thm:oneleg}. If
$c=2$ let $\lambda=\lambda_{\phi'(\partial_1), \phi'(\partial_2)}$ be as in theorem
\ref{thm:twoleg}. Theorems \ref{thm:oneleg} and \ref{thm:twoleg} imply that for $\phi'_t
\in R(E_t)$
\begin{eqnarray*}
|\pi^{-1}(R^{-1}(\phi') \cap E_t)| \sim \lambda(R^{-1}(\phi') \cap E_t) =
\lambda(E_{t,\phi'}).
\end{eqnarray*}

Using the induction hypothesis, this implies
\begin{eqnarray*}
|\pi^{-1}(E_t)| &=& \sum_{[\phi'] \in R(E_t)} \, |\pi^{-1}(R^{-1}(\phi') \cap E_t)| \sim
\int_{R(E_t)} \lambda(E_{t,\phi'}) \, d\nu'(\phi').
\end{eqnarray*}
A routine calculation shows that this equals $\nu(E_t)$.

\end{proof}

\begin{proof}(of Theorem \ref{thm:2d})
It needs only be noted that the sequence $\{N_\epsilon [\phi_t]\}_t$ satisfies the
hypothesis of theorem \ref{thm:3d}. This is easy and we leave it to the reader.
\end{proof}

%\section{Notes}

%People to send the paper to: Bill Goldman, Yiannis Petridis, Alex Eskin, Chris Leininger,
%Chris Connell, Chris Judge, David Fisher, the guy in Brazil, Dolgopyat, Forni, Motohico.

%Reminders: The formulas generalize the classical pgt.

\end{document}